\documentclass[11pt]{article}

\usepackage{hyperref}

\usepackage{amsmath,amsfonts,amssymb,amsthm}
\usepackage{mathrsfs}
\usepackage{enumitem}
\usepackage{graphicx}
\usepackage{color}

\newtheorem{theorem}{Theorem}[section]

\newtheorem{proposition}[theorem]{Proposition}

\theoremstyle{remark}
\newtheorem{remark}[theorem]{Remark}
\renewenvironment{proof}[1][Proof]{\noindent{\itshape {#1.} } }{$\Box$ 
\medskip} 

\numberwithin{equation}{section}
\newcommand{\R}{\mathbb{R}}
\newcommand{\Z}{\mathbb{Z}}

\newcommand{\E}{\mathbb{E}}
\newcommand{\Cc}{\mathcal{C}}
\newcommand{\eps}{\varepsilon}
\newcommand{\sinc}{\mathrm{sinc}}

\newcommand{\sigref}{\sigma_{\rm r}}
\newcommand{\xref}{x_{\rm r}}
\newcommand{\cref}{c_{\rm r}}
\newcommand{\xs}{x_{\mathrm{s}}}
\newcommand{\zr}{z_{\mathrm{r}}}
\newcommand{\sA}{\mathcal{A}}

\title{Passive Array Correlation-Based Imaging in a Random Waveguide}
\author{Habib Ammari\thanks{D\'epartement de Math\'ematiques et Applications, Ecole Normale Sup\'erieure, 
45 Rue d'Ulm, 75230 Paris Cedex 05, France. Email: {\tt habib.ammari@ens.fr}} \and Josselin Garnier\thanks{Laboratoire de
Probabilit\'es et Mod\`eles Al\'eatoires \& Laboratoire
Jacques-Louis Lions, Universit\'e Paris VII, 75205 Paris Cedex 13,
France. Email: {\tt garnier@math.univ-paris-diderot.fr}}  \and 
Wenjia Jing\thanks{D\'epartement de Math\'ematiques et Applications, Ecole Normale Sup\'erieure, 
45 Rue d'Ulm, 75230 Paris Cedex 05, France. Email: {\tt wjing@dma.ens.fr}}}
\begin{document}
\maketitle
\begin{abstract}
We consider reflector imaging in a weakly random waveguide. We address the situation in which 
the source is farther from the reflector to be imaged than the energy equipartition distance, but the 
receiver array is closer to the reflector to be imaged than the energy equipartition distance.
As a consequence, the reflector is illuminated by a partially coherent field and the signals recorded by the receiver array are noisy.
This paper shows that migration of the recorded signals cannot give a good image, but an appropriate
migration of the cross correlations of the recorded signals can give a very good image.
The resolution and stability analysis of this original functional shows that the reflector can be localized with
an accuracy of the order of the wavelength even when the receiver array has small aperture,
and that broadband sources are necessary to ensure statistical stability, whatever the aperture of the array.
\end{abstract}

%%%%%%%%%%%%%%%
%%%%%%%%%%%%%%%
\section{Introduction}
\label{sec:intro}
Sensor array imaging in a scattering medium is limited because 
coherent signals recorded at the source-receiver array and coming
from a reflector to be imaged are dominated
by incoherent signals coming from multiple scattering by the medium.
For instance, in a randomly perturbed waveguide, it is known that the field becomes completely incoherent
when the propagation distance becomes larger than the equipartition distance, which corresponds to the distance beyond which the
source energy has been shared equally among all the propagating modes \cite[Chapter 20]{FGPS}.
As we will see, if the distance between the source-receiver array and the reflector is larger than the equipartition distance, then
classical migration of the signals recorded at the array cannot give a good image.

Sources can be expensive or difficult to implement 
but receivers can be cheap and easy to implement, so an imaging problem in which there are a few sources
(all of them being far from the reflector) and many receivers (some of them being close to the reflector)
is of theoretical and practical  interest.
If there is a unique source far from the reflector (farther than the equipartition distance)
and if the receiver array is close to the reflector (closer than the equipartition distance),
then  classical migration of the recorded signals fails again. This was shown is various contexts and we will show it again in the waveguide geometry.
However, in such a situation, another kind of migration can be used: from the work devoted to coherent interferometry imaging
\cite{borcea11a,borcea11b,borcea05,borcea06a,borcea06b} and ambient noise imaging \cite{curtis06,cross,sabra05,schuster,wap10},
it is known that migration of cross correlations of noisy signals can be more stable than migration of the signals themselves.
The migration of cross correlations of noisy signals recorded by auxiliary passive arrays was proposed by \cite{bakulin} in geophysical contexts and 
analyzed recently in randomly scattering open media in \cite{virtual},
and we would like to address the same problem in the waveguide geometry.
Indeed the number of propagating modes is finite in the waveguide geometry so that the statistical behavior of partially coherent fields
in random waveguides is very different from the open medium case \cite{GP-SIAM07,FGPS}.
In our paper, we show that, if a receiver array can be placed close to the reflector
to be imaged, then the cross correlations of the
incoherent signals on this array can  be used to image the reflector.
We will give a detailed resolution and stability analysis. We will show that
the statistical stability requires a broadband source 
 and that good resolution and stability properties do not require the receiver array to span the whole cross section of the waveguide,
which is an effect specific to the waveguide geometry.

The  paper is organized as follows. In section \ref{sec:prelim}, we review the mathematical background of the imaging problem in a random waveguide.
In section \ref{sec:imagclas} we describe and analyze the classical migration functional using the recorded signals and show that 
it cannot give a good image when the propagation distance is beyond the energy equipartition distance.
In section \ref{sec:imagf}, we introduce the correlation-based imaging functional; 
it has two versions which correspond to the time-harmonic case and the broadband case. 
In section \ref{sec:resol}, we analyze the resolution of the proposed imaging functionals. 
Detailed analyses are provided for full aperture and limited aperture arrays. These results are based on the statistical average of the imaging functionals. The variances of these functionals are very important as well because they determine the statistical stability of the imaging functionals. 
In section \ref{sec:stabi}, we study the variances of the imaging functionals. Some concluding remarks are listed at the end of the paper.
 
%%%%%%%%%%%%%%%
%%%%%%%%%%%%%%%
\section{Mathematical Formulation of the Imaging Problem}
\label{sec:prelim}

%%%%%%%%%%%%%%%
\subsection{The ideal waveguide}
\label{sec:acous}

We consider linear scalar (acoustic) waves propagating in a two-dimensional space. The governing equation is
\begin{equation}
\Delta p(t,\mathbf{x}) - \frac{1}{c^2_0} \frac{\partial^2 p}{\partial t^2} (t,\mathbf{x}) = F (t,\mathbf{x}).
\label{eq:pwe}
\end{equation}
Here $p$ is the scalar field (acoustic pressure);
$c_0$ is the speed of propagation in the medium (sound speed);
$F(t,\mathbf{x})$ models the forcing term. 
We consider a waveguide geometry, and we decompose the spatial variable $\mathbf{x}$ as $(x,z)$. That is, $z \in \R$ is along the axis of the waveguide while $x \in {\mathcal D}$ denotes the transverse coordinate, and ${\mathcal D}= (0,a)$ is the transverse section of the waveguide. 
We assume that the forcing term is localized in the plane $z=0$:
\begin{equation}
F(t,\mathbf{x}) = f(t) \delta(\mathbf{x} - \mathbf{x}_{\mathrm{source}}) ,
\label{eq:Fdef}
\end{equation}
where $\mathbf{x}_{\mathrm{source}} = (\xs,0)$ for some $\xs \in \mathcal{D}$.
We assume that the medium is quiescent before the pulse emission, that is
\begin{equation}
p(t,\mathbf{x}) = 0, \quad t\ll 0.
\end{equation}
We consider Dirichlet boundary conditions at the boundary of the waveguide:
\begin{equation}
p(t,\mathbf{x}) = 0, \quad \mathbf{x} \in \partial \mathcal{D} \times \R =\{0,a\} \times \R.
\end{equation}

Using the Fourier method, the scalar field can be written as a superposition of waveguide modes.
A waveguide mode is a time-harmonic wave of the form $\hat{p}(\omega,\mathbf{x}) e^{-i\omega t}$ with frequency $\omega$, where $\hat{p}$ satisfies the time-harmonic form of the wave equation \eqref{eq:pwe} without a source term:
\begin{equation}
\partial_z^2 \hat{p}(\omega,x,z) + \Delta_{\perp} \hat{p}(\omega,x,z) +  
 k^2(\omega) \hat{p}(\omega,x,z) = 0.
\label{eq:pbarwe}
\end{equation}
Here, $\Delta_\perp = \partial_x^2$ is the transverse Laplace operator in the transverse section $\mathcal{D}$ with Dirichlet boundary conditions; $k(\omega) = \omega/c_0$ is the homogeneous wavenumber. 
Consequently, \eqref{eq:pbarwe} can be solved using the eigenmodes of $\Delta_\perp$, that is, using the orthonormal basis $\{\phi_j(x)\}_{j=1,2,\cdots}$ of $L^2(\mathcal{D})$ given by
\begin{equation}
-\Delta_\perp \phi_j(x) = \lambda_j \phi_j(x), \quad x \in \mathcal{D}.
\end{equation}
The eigenvalues are simple, satisfying $0<\lambda_1 < \lambda_2 < \cdots$.
The eigenvalues and eigenvectors are given by
\begin{equation}
\phi_j (x)= \frac{\sqrt{2}}{\sqrt{a}} \sin(\frac{\pi j x}{a}), \quad \lambda_j = \frac{\pi^2 j^2}{a^2} .
\end{equation}
Using the method of separation of variables on \eqref{eq:pbarwe}, we see that the waveguide mode $\hat{p}(\omega,x,z)$ can be further written as superposition of $\hat{p}_j(\omega,x,z) = \phi_j(x) e^{\pm i\beta_j (\omega) z}$, where
\begin{equation}
\beta_j^2 (\omega) = k^2(\omega) - \lambda_j.
\end{equation}
For a given frequency $\omega$, there exists a unique integer $N(\omega)$ such that $\lambda_{N(\omega)} \le k^2(\omega) < \lambda_{N(\omega)+1}$:
\begin{equation}
N(\omega) = \Big\lfloor\frac{\omega a}{\pi c_0} \Big\rfloor.
\end{equation}
Here and in the sequel, $\lfloor b \rfloor$ means the integer part of a real number $b$.
 The modes $\{\hat{p}_j(\omega,x,z) = \phi_j(x) e^{\pm i\beta_j(\omega) z}\}_{j=1,\ldots,N(\omega)}$ are propagating waveguide modes and $\{\beta_j(\omega)\}_{j=1,\ldots,N(\omega)}$ are called the modal wavenumbers. On the other hand, $\{\hat{p}_j(\omega,x,z) = \phi_j(x) e^{\pm |\beta_j(\omega)|z}\}_{j>N(\omega)}$ are evanescent modes because they decay as $z$ goes to $\mp \infty$.

%%%%%%%%%%%%%%%
\subsection{The randomly perturbed waveguide}
\label{sec:rands}%
From now on we assume that the waveguide is randomly perturbed and the scalar field satisfies the perturbed wave equation
\begin{equation}
\Delta p(t,\mathbf{x})  - \frac{1}{c^2(\mathbf{x})} \frac{\partial^2 p}{\partial t^2} (t,\mathbf{x}) = F (t,\mathbf{x}),
\label{eq:pwepert}
\end{equation}
where $c(\mathbf{x})$ is the randomly heterogeneous speed of propagation of the medium.
We consider the case where the typical amplitude of the fluctuations of the speed of propagation is small, which we call the weakly random regime. 
When the correlation length of the fluctuations is of the same order as the typical wavelength 
the interactions between the waves and the random medium become nontrivial. 
Due to the small amplitude of the fluctuations, however, the effect of the interaction becomes important only after a long propagating distance.

More exactly we assume that a randomly heterogeneous section in $z \in [0,\widetilde{L}_0]$ is sandwiched in between two homogeneous waveguides:
The speed of propagation is of the form
\begin{equation}
\frac{1}{c^2(x,z)} =\left\{
\begin{array}{ll}
 \frac{1}{c_0^2} (1+\eps \nu(x,z) ) &\mbox{ if }(x ,z) \in  [0,a] \times [0,\widetilde{L}_0],\\
 \frac{1}{c_0^2} &  \mbox{otherwise}
\end{array}
\right. 
\label{eq:rhoKdef}
\end{equation}
Here, $\nu$ is a mean-zero, stationary and ergodic random processes with respect to the axis coordinate $z$. 
It is assumed to satisfy strong mixing conditions in $z$.
The relative amplitude of the fluctuations of the speed of propagation is denoted by $\eps$.
We assume that the correlation length of the random perturbation is of the same order as the typical wavelength  
$\lambda(\omega_0)= 2\pi c_0 /\omega_0 = 2\pi / k(\omega_0)$, for $\omega_0$ the central frequency of the source.
We assume that the propagation distance $\widetilde{L}_0$ is much larger than the typical
wavelength. We will see that the interesting regime is when the ratio $\lambda/\widetilde{L}_0$ is of order $\eps^2$,
so we introduce the normalized propagation distance $L_0$:
$$
\widetilde{L}_0=\frac{L_0}{\eps^2} .
$$
In this regime the cumulative effects of the interaction of the scalar wave with the small fluctuations of the speed of propagation
become of order one.

For a fixed frequency $\omega$, 
 the Fourier transformed scalar field $\hat{p}(\omega,x,z)$ defined by
 $$
 \hat{p}(\omega,x,z) = \int {p}(t,x,z) e^{i \omega t} dt
 $$
 satisfies the equation
\begin{equation}
\partial_z^2 \hat{p}(\omega,x,z) + \Delta_{\perp} \hat{p}(\omega,x,z) + k^2(\omega)[1+\eps \nu(x,z)] \hat{p}(\omega,x,z) = \hat{f}(\omega) \delta(x-x_{\rm s}) \delta(z) .
\label{eq:pbarwer}
\end{equation}
To solve this equation, we make the following two simplifications that are justified in \cite[Chapter 20]{FGPS} or \cite{GP-SIAM07}.

\emph{Ignoring the evanescent modes}. 
First, we only consider the propagating modes:
\begin{equation}
\hat{p}(\omega,x,z) = \sum_{j=1}^{N(\omega)} \phi_j(x) 
\hat{p}_j(\omega,z).
\label{eq:ptrunc}
\end{equation}
This is valid because we are mainly concerned with the scalar field for $z \gg 1$ and the evanescent modes decay exponentially fast. Furthermore,
we parameterize the complex mode amplitude $\hat{p}_j(\omega,z)$ by the amplitudes of its right- and left-going components. Let $\hat{a}_j(\omega,z)$ and $\hat{b}_j(\omega,z)$ be the amplitudes of these components, defined by
\begin{equation}
\hat{p}_j = \frac{1}{\sqrt{\beta_j}} \left(\hat{a}_j e^{i\beta_j z} + 
\hat{b}_j e^{-i\beta_j z} \right),
\quad\quad
\frac{d\hat{p}_j}{dz} = i\sqrt{\beta_j} \left(\hat{a}_j e^{i\beta_j z} - 
\hat{b}_j e^{-i\beta_j z} \right).
\label{eq:leftright}
\end{equation}
Using these representations, one obtains a system of ordinary differential equations (ODEs) for $\{\hat{a}_j,\hat{b}_j\}$
 \cite[Section 20.2.4]{FGPS} or \cite[Section 3.1]{GP-SIAM07}. The coefficients of the system depend on the integrated quantities of the form
\begin{equation}
C_{jl}(z) = \int_{\mathcal{D}} \phi_j(x) \phi_l(x) \nu (x,z) dx.
\end{equation}
This system of ODEs is closed by the boundary conditions at $z=0$ where the source $F$ is imposed and at $z=\widetilde{L}_0=L_0/\eps^2$ where there is no left-going component.

\emph{Forward scattering approximation}. Second, we neglect the left-going (backward) propagating mode, assuming that they do not interact with the right-going ones. This is valid in the limit $\eps \to 0$ when the second-order moments of 
$\nu$ satisfy certain conditions  \cite[Section 20.2.6]{FGPS} or \cite[Section 3.3]{GP-SIAM07}. In this case, the rescaled amplitude 
$$
\hat{a}^\eps_j (z)= \hat{a}_j(\omega,z/\eps^2)
$$
 of the right-going wave satisfies
\begin{equation}
\frac{d \hat{\bf a}^\eps}{dz} = \frac{1}{\eps} \mathbf{H}^{(a)}_\omega \left(\frac{z}{\eps^2}\right) \hat{\bf a}^\eps,
\label{eq:aeq}
\end{equation}
where $\hat{\bf a}^\eps$ denotes the $N(\omega)$-dimensional vector $(\hat{a}^\eps_1,\ldots,\hat{a}^\eps_{N(\omega)})'$ and $\mathbf{H}^{(a)}_\omega$ is a $N(\omega)\times N(\omega)$ complex matrix with components
\begin{equation}
H^{(a)}_{\omega,jl} = \frac{ik^2}{2} \frac{C_{jl}(z)}{\sqrt{\beta_j \beta_l}} e^{i(\beta_l-\beta_j)z}.
\label{eq:Haomega}
\end{equation}
Define the propagator matrix $\mathbf{T}^\eps(\omega,z,z_0)$ to be the fundamental solution of the system \eqref{eq:aeq}, {\itshape i.e.},
\begin{equation}
\frac{d \mathbf{T}^\eps}{dz} (\omega,z,z_0) = \frac{1}{\eps} \mathbf{H}^{(a)}_\omega \left(\frac{z}{\eps^2}\right) \mathbf{T}^\eps(\omega,z,z_0),
\label{eq:Teq}
\end{equation}
with $\mathbf{T}^\eps(\omega,z=z_0,z_0) = \mathbf{I}$. Then $\hat{\bf a}^\eps(\omega,z) = \mathbf{T}^\eps(\omega,z,0) \hat{\bf a}^\eps(\omega,0)$, where the initial amplitude $\hat{\bf a}^\eps(\omega,0)$ is determined by the source $F$. In fact, 
integrating (\ref{eq:pbarwer}) across the plane $z=0$ and using  (\ref{eq:pbarwer},\ref{eq:leftright}), we find that
\begin{equation*}
\hat{a}_l(\omega,0) = \frac{1}{2 i  \sqrt{\beta_l(\omega)}} \hat{f}(\omega) \phi_l(\xs), \quad l=1,\ldots,N(\omega).
\end{equation*}

Consider an array of receivers located in the plane 
$z=\widetilde{L}$  of the random waveguide section,
where
\begin{equation}
\widetilde{L}=\frac{L}{\eps^2},
\end{equation}
and $0 < L < L_0$. 
Let $T^\eps_{jl}(\omega)$ 
 be the $jl$-entry of the propagator matrix $\mathbf{T}^\eps(\omega,L,0)$. 
It is the rate of conversion  of the initial $l$-mode into the $j$-mode
in the plane $z=\widetilde{L}=L/\eps^2$ of the random waveguide section. In particular, we have
\begin{equation}
\hat{a}_j(\omega, \widetilde{L}) = \hat{a}^\eps_j (\omega,L) = \sum_{l=1}^{N(\omega)} T^\eps_{jl}(\omega) \hat{a}^\eps_l(\omega,0) = \sum_{l=1}^{N(\omega)} \frac{1}{2 i \sqrt{\beta_l(\omega)}}
T^\eps_{jl}(\omega) \hat{f}(\omega) \phi_l(\xs).
\label{eq:ahatj}
\end{equation}
Repeating the argument above, we see that the field beyond $\widetilde{L}$, that is 
$\hat{\bf a}(\omega,z) = \{\hat{a}_j(\omega,z)\}_{j=1}^{N(\omega)}$ for $z> \widetilde{L}=L/\eps^2$, are related to $\hat{\bf a}(\omega,\widetilde{L}) = \{\hat{a}_j(\omega,\widetilde{L})\}_{j=1}^{N(\omega)}$ as follows
$$
\hat{a}_j(\omega,z) = \sum_{l=1}^{N(\omega)}  {T}^\eps_{jl}(\omega,\eps^2 z, L) \hat{a}^\eps_l(\omega,L).
$$
Since the random waveguide is stationary, $\{{T}^\eps_{jl}(\omega,z,z_0)\}_{j,l=1}^{N(\omega)}$ has the same distribution as $\{T^\eps_{jl}(\omega,z-z_0,0)\}_{j,l=1}^{N(\omega)}$. Therefore, we can apply \eqref{eq:prop:T1} and \eqref{eq:Tsqlim} in Proposition \ref{prop:momentT} and conclude that ${T}^\eps_{jl}(\omega,\eps^2 z, L) \approx \delta_{jl}$ in probability provided that $z-\tilde{L} \ll \eps^{-2}$. This is equivalent to say
$$
\hat{a}_j(\omega, z) \approx \hat{a}_j(\omega, \widetilde{L}), \quad\quad \text{for } 0\le z - \widetilde{L} \ll \eps^{-2}.
$$

Using this approximation and the expressions \eqref{eq:ptrunc} and \eqref{eq:leftright}, we can write the scalar field at $z > \widetilde{L}$ with $z-\widetilde{L} \ll \eps^{-2}$ as:
\begin{equation}
\hat{p}   (\omega,x,z)  = \sum_{j,l=1}^{N(\omega)}
\frac{1}{2 i  \sqrt{\beta_l(\omega)} \sqrt{\beta_j(\omega)}} 
T^\eps_{jl}(\omega) \hat{f}(\omega) \phi_j(x) \phi_l(\xs) e^{i\beta_j z}.
\label{eq:ppfs0}
\end{equation}
%[END OF REVISION]

\begin{figure}
\vspace*{-1.6cm}
\begin{center}
\includegraphics[width=9cm]{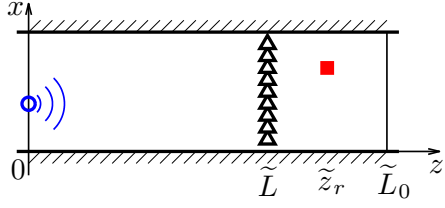}
\end{center}
\vspace*{-2.6cm}
\caption{Schematic of the imaging problem. A point source (circle) in the plane $z=0$ 
emits a short pulse that propagates through the random waveguide.
The target (square) in the plane $z=\widetilde{z}_{\rm r}$ is a reflector. The receiver array (triangles) 
 in the plane $z=\widetilde{L}$ records the signals. }
\label{fig0}
\end{figure}

%%%%%%%%%%%%
\subsection{Modeling the point reflector}
\label{subsec:reflector}%
In the imaging problem to be investigated (see Figure \ref{fig0}), the goal is to locate a point reflector
centered at $\mathbf{x}_{\rm r} = (x_{\rm r},\widetilde{z}_{\rm r})$ from signals recorded at the receiver array in the plane 
$z=\widetilde{L}=L/\eps^2$. 
The reflector is supposed to be at a relatively small distance (compared to $\eps^{-2}$), that is to say
%[REVISION: ADDED Z_{\RM_R} \LL EPS^{-2}]
\begin{equation}
\widetilde{z}_{\rm r}=\frac{L}{\eps^2}+z_{\rm r}, \quad\quad 1\ll z_{\rm r} \ll \eps^{-2}.
\end{equation}
Note that we also assume that $z_{\rm r} \gg 1$, i.e. the distance between the reflector and the receiver array is much larger than one (the order of magnitude of the wavelength), to ignore the evanescent modes emitted by the reflector.
%[END OF REVISION]
The reflector can be modeled as a local change in the density and/or the bulk modulus of the medium, so that the sound speed is locally modified as
\begin{equation}
\frac{1}{c^2(x,z)} = \frac{1}{c_0^2} \big(1 + \eps \nu(x,z) \big) +  \frac{1}{\cref^2} \mathbf{1}_{\Omega_{\rm r}}(x,z) ,
\label{eq:mref}
\end{equation}
where ${\Omega_{\rm r}}$ is a small domain around $\mathbf{x}_{\rm r}:= (\xref,\widetilde{z}_{\rm r})$ which represents the center of the reflector; $\cref$ is a parameter characterizing the contrast of the reflector. With this modification, the right-hand side of \eqref{eq:pbarwe} should have an additional term $-(\omega/\cref)^{2}\mathbf{1}_{\Omega_{\rm r}} \hat{p}(\omega,x,z)$. We assume that the diameter of the scattering region ${\Omega_{\rm r}}$ is small compared to the typical wavelength and that the velocity contrast is such that $\sigref: = \cref^{-2}|{\Omega_{\rm r}}|$ satisfies 
 $\sigref \ll 1$. Then we can model the scattering region by a point reflector
\begin{equation*}
\frac{1}{\cref^2} \mathbf{1}_{\Omega_{\rm r}}(x,z) \approx \sigref \delta(\mathbf{x} - \mathbf{x}_{\rm r}).
\end{equation*}

\emph{Born approximation}. The above setting allows us to solve the scalar field 
with the presence of the point reflector using the Born approximation for the reflector. Given a fixed frequency $\omega$, we have
\begin{equation}
\label{eq:phatdec}
\hat{p} (\omega,x,z) \approx \hat{p} _{\rm p} (\omega,x,z) + \hat{p}_{\rm s} (\omega,x,z).
\end{equation} 
Here, $\hat{p}_{\rm p} $ is the primary field induced by the source ${F}$ propagating through the random waveguide 
and computed in the previous section (Eq.~(\ref{eq:ppfs0})), 
and $\hat{p}_{\rm s} $ is the secondary field, that is the first-order scattered field due to the additional source $-\omega^2 \sigref \delta(\mathbf{x}-\mathbf{x}_{\rm r}) \hat{p}_{\rm p} $ at the reflector:
\begin{equation}
\Delta \hat{p}_{\rm s}  (\omega, x, z) + k^2(\omega) \big[1 +\eps \nu(x,z)\big] \hat{p}_{\rm s}  
(\omega, x, z) = -\omega^2\sigref \delta(\mathbf{x} - \mathbf{x}_{\rm r}) \hat{p}_{\rm p} (\omega, x_{\rm r},\widetilde{z}_{\rm r}).
\label{eq:pbars}
\end{equation}
 Note that in the Born approximation one replaces the full wave field at the reflector by the primary field in the right-hand side of (\ref{eq:pbars}).

The primary field is solved exactly as in the previous section. Summarizing the results there, one obtains 
that for $z \ge \widetilde{L}$ and $z -\widetilde{L} \ll \eps^{-2}$,
\begin{equation}
\hat{p}_{\rm p}  (\omega,x,z) 
=
 \sum_{j,l=1}^{N(\omega)}
\frac{\hat{f}(\omega) }{2 i  \sqrt{\beta_l(\omega)} \sqrt{\beta_j(\omega)}} 
T^\eps_{jl}(\omega)  \phi_j(x) \phi_l(\xs) e^{i\beta_j z}.
\label{eq:ppfs}
\end{equation}

The secondary field satisfies (\ref{eq:pbars}).
Again, we solve this equation using the orthonormal basis $\{\phi_j(x)\}_{j=1,\ldots,N(\omega)}$ and we ignore the evanescent modes. 
Since the reflector is within a distance smaller than $\eps^{-2}$ from the receiver array, for $\widetilde{L} < z < \widetilde{z}_{\rm r}$, the propagator matrix from $\widetilde{z}_{\rm r}$ to $z$ can be approximated by the  identity matrix in probability, and we only need to decompose the secondary source at the reflector into waveguide modes.
Using \eqref{eq:pbars}, the decomposition \eqref{eq:leftright} and the fact that there is no left-going wave from $z>\widetilde{z}_{\rm r}$, we find that
\begin{equation}
\hat{b}_{{\rm s}j}(\omega, z) = \frac{i\omega^2 \sigref}{2\sqrt{\beta_j}}  \phi_j(\xref) \hat{p}_{\rm p} (\omega, \xref,\widetilde{z}_{\rm r}).
\label{eq:secondb}
\end{equation}
We note that there is no right-going secondary wave because we do not consider back-scattering of the left-going secondary wave near the receivers. Finally, recall the expression of the primary field at the reflector \eqref{eq:ppfs}, we obtain
 for $z\in [ \widetilde{L} , \widetilde{z}_{\rm r})$ that
\begin{eqnarray}
\nonumber
\hat{p}_{\rm s} (\omega, x, z) &=& \sum_{j=1}^{N(\omega)} 
\phi_j(x) \frac{1}{\sqrt{\beta_j}} \hat{b}_{{\rm s}j}(\omega,z) e^{-i\beta_j(z-\widetilde{z}_{\rm r})}\\
&= &\sum_{j,l,m=1}^{N(\omega)} \frac{\omega^2 \sigref \hat{f}(\omega) }{4\beta_j 
 \sqrt{\beta_m} \sqrt{\beta_l}} T^\eps_{lm}(\omega)
\phi_j(x) \phi_j(\xref) \phi_l(\xref) \phi_m(\xs) e^{-i\beta_j(z-\widetilde{z}_{\rm r})} e^{i\beta_l \widetilde{z}_{\rm r}}. \hspace*{0.2in}
\label{eq:psfs}
\end{eqnarray}
%Since $\sigref \ll 1$ is a very small number, the secondary field field is dominated by the primary field.

%%%%%%%%%%%%%%%
%%%%%%%%%%%%%%%

\section{Migration-Based Imaging Functional}
\label{sec:imagclas}
In this section, we introduce the classical imaging functional to localize the point reflector 
using the scalar (pressure) field recorded at the receiver array at $z=\widetilde{L}$. 
This imaging functional is based on the migration of the array data to a search point 
$(x^{\rm S},\widetilde{z}^{\rm S})$. 
Our goal is to show that classical Kirchhoff migration functional does not give a good image 
when the medium between the source at $z=0$ 
and the receiver array at $z=\widetilde{L}$ is scattering.

The data of scalar (pressure) field recorded by the receivers are 
$$
\{p(t,x,\widetilde{L}) ~|~ t\in \R, x \in {\mathcal D}\} .
$$
Note that we consider in this section the full aperture case:
the receivers span the whole cross section of the waveguide and they record data at all time.
We consider the frequency- and mode-dependent data
$$
\hat{p}_j(\omega,z=\widetilde{L}) =\int \int p (t,x,z=\widetilde{L})\phi_j(x) dx e^{ i\omega t} dt  .
$$
According to the analysis carried out Section \ref{subsec:reflector},
it can be decomposed as
$$
\hat{p}_j(\omega,z=\widetilde{L})= \hat{p}_{{\rm p}j}(\omega,z=\widetilde{L}) +
\hat{p}_{{\rm s}j}(\omega,z=\widetilde{L})   .
$$
From (\ref{eq:ppfs}) and (\ref{eq:psfs}) the primary and secondary contributions are
\begin{eqnarray}
\hat{p}_{{\rm p}j}(\omega,z=\widetilde{L}) 
&=&
 \frac{\hat{f}(\omega)}{2 i   \sqrt{\beta_j(\omega)}}e^{i\beta_j \frac{L}{\eps^2}}
\sum_{l=1}^N  \frac{T^\eps_{jl}(\omega)}{ \sqrt{\beta_l}} \phi_l (\xs)  ,
\\
\hat{p}_{{\rm s}j}(\omega,z=\widetilde{L}) 
&=&  \frac{1}{\beta_j(\omega)}  \phi_j(\xref) e^{i\beta_j  {z}_{\rm r}} q(\omega, x_{\rm r},\widetilde{z}_{\rm r}) ,
\label{eq:psj}
\end{eqnarray}
with 
$$
q(\omega, x_{\rm r},\widetilde{z}_{\rm r}) = \sum_{l,m=1}^N
 \frac{\omega^2\hat{f}(\omega)  \sigref}{4  \sqrt{\beta_m(\omega)} \sqrt{\beta_l(\omega)}} T^\eps_{lm}(\omega)
\phi_l(\xref) \phi_m(\xs)  e^{i\beta_l (\omega)\widetilde{z}_{\rm r}},
$$
which can be interpreted as an illumination of the reflector.
The secondary contribution $\hat{p}_{{\rm s}j}$ contains the information about the reflector,
and its form (\ref{eq:psj}) motivates the definition of the Kirchhoff migration imaging functional:
\begin{equation}
\mathcal{I}_{\mathrm{KM}}(x^{\rm S}, z^{\rm S}) = \frac{1}{2\pi} \int 
\frac{1}{N(\omega)} \sum_{j=1}^{N(\omega)} \beta_j (\omega)\phi_j (x^{\rm S}) e^{-i \beta_j(\omega) z^{\rm S}}
 \hat{p}_j(\omega,z=\widetilde{L}) d\omega ,
\end{equation}
where the search point is $(x^{\rm S},\widetilde{z}^{\rm S})$ with $\widetilde{z}^{\rm S}=\widetilde{L}+z^{\rm S}$.

A simple case is when the source term is time-harmonic, {\itshape i.e.} $F(t,\mathbf{x})=\delta(\mathbf{x}-\mathbf{x}_{\mathrm{source}}) f(t) $
with $f(t) = e^{-i\omega_0 t}$ and
\begin{equation*}
\hat{f}(\omega) = 2\pi \delta(\omega-\omega_0).
\end{equation*}
Then the data set is reduced to $\{ \hat{p}_j(\omega_0,z=\widetilde{L}) , j =1,\ldots,N(\omega_0)\}$ and the Kirchhoff migration functional has the form
\begin{equation}
\mathcal{I}_{\mathrm{KM}}(x^{\rm S}, z^{\rm S}) =  
\frac{1}{N(\omega_0)} \sum_{j=1}^{N(\omega_0)} \beta_j (\omega_0)\phi_j (x^{\rm S}) e^{-i \beta_j(\omega_0) z^{\rm S}} \hat{p}_j(\omega_0,z=\widetilde{L})   .
\end{equation}
We need to compute the mean of the imaging functional   in order to characterize its resolution properties 
and its variance  in order to characterize its stability properties.
These statistical moments depend on the moments of the propagator matrix which were studied in \cite[Propositions 20.6 and 20.8]{FGPS}
or \cite[Propositions 6.1 and 6.3]{GP-SIAM07}.

\begin{proposition}
\label{prop:momentT} 
The first-order moments of the transmission coefficients have limits as $\eps \to 
0$, which are given by
\begin{align}
\E [T^\eps_{jl}(\omega) ] 
&\xrightarrow[]{\eps \to 0} 0, & &\text{ if } j\ne
l,\\
\E [T^\eps_{jj}(\omega)  ] 
&\xrightarrow[]{\eps \to 0} e^{ - D_{j}(\omega)L}, & &  \text{ otherwise.}
\end{align}
The second-order moments of the
transmission coefficients  have limits as $\eps \to 
0$, which are given by
\begin{align}
\E [T^\eps_{jj}(\omega) \overline{T^\eps_{ll}}(\omega)] 
&\xrightarrow[]{\eps \to 0} e^{- Q_{jl}(\omega)L}, & & \text{ if } j\ne 
l,\\
\E [T^\eps_{jl}(\omega) \overline{T^\eps_{jl}}(\omega)] 
&\xrightarrow[]{\eps \to 0} \mathcal{T}^{(l)}_j(\omega,L), & &\label{eq:prop:T1} \\
\E [T^\eps_{jl}(\omega) \overline{T^\eps_{mn}}(\omega)] 
&\xrightarrow[]{\eps \to 0} 0, & & \text{ otherwise.}
\end{align}
The functions $\mathcal{T}_j^{(l)}(\omega, z)$ are
the solutions of the system of linear equations
\begin{equation}
\frac{d \mathcal{T}_j^{(l)}}{dz}  = \sum_{n\ne j} 
\Gamma^{(c)}_{jn}(\omega) \left(\mathcal{T}^{(l)}_n - \mathcal{T}^{(l)}_j 
\right), \quad \quad \mathcal{T}^{(l)}_j(\omega, z = 0) = \delta_{jl}.
\label{eq:Tsqlim}
\end{equation}
The positive coefficients $D_{j}$ and $Q_{jl}$ and the matrix $\Gamma^{(c)}_{jn}$ depend on the 
correlation function of the random process $\nu$. Furthermore, we have
\begin{equation}
\sup_{j,l} \left\lvert \E [T^\eps_{jl}(\omega) ] 
\right\rvert \le C e^{-L/L_{\mathrm{equip}}}, \quad \quad
\sup_{j,l} \left\lvert \mathcal{T}_j^{(l)} (\omega,L) - \frac{1}{N} 
\right\rvert \le C e^{-L/L_{\mathrm{equip}}},
\label{eq:sTjllim}
\end{equation}
where $L_{\mathrm{equip}}$ is the equipartition distance for the mean 
mode powers introduced at the end of Section 20.3.3 in \cite{FGPS}
(or at the end of Section 4.2 in \cite{GP-SIAM07}).
\end{proposition}
The results on the first-order moments describe how the wave loses its coherence as it propagates in the random waveguide.
The results on the second-order moments describe how the wave energy becomes equipartitioned on the waveguide modes.
 
When $L$ is larger than the energy equipartition length 
$L_{\mathrm{equip}}$, then the first-order moments of the transmission coefficients are vanishing.
Based on this observation, we  have 
\begin{equation}
\E \big[ \mathcal{I}_{\mathrm{KM}}(x^{\rm S}, z^{\rm S}) \big] \approx 0 .
\end{equation}

It turns out that the fluctuations of the imaging functional are much larger than its mean.
This can be seen by studying the standard deviation of the imaging functional.
When $L$ is larger than the energy equipartition length 
$L_{\mathrm{equip}}$, then the second-order moments 
of the transmission coefficients are vanishing
except $\E [|T_{jl}^\eps|^2]$ which converge to $1/N$.
Based on this observation, the second-order moment of the imaging functional for a time-harmonic source is:
\begin{eqnarray}
\nonumber
\E \big[ | \mathcal{I}_{\mathrm{KM}}(x^{\rm S}, z^{\rm S}) |^2 \big] 
&=& \frac{|\hat{f}(\omega_0)|^2 \Phi_{-1}(x_{\rm s})}{N} \Big[
\frac{1}{4}   \Phi_{1}(x^{\rm S})
+
 \Big( \frac{\omega^2 \sigma_{\rm r} N}{4}\Big)^2 \Phi_{-1} (x_{\rm r})   
|\Psi (x^{\rm S},z^{\rm S};\xref,\zr)|^2 \\
&& +
 \Big( \frac{\omega^2 \sigma_{\rm r} N}{4}\Big)    
\Im m \big( {\Psi (x^{\rm S},-z^{\rm S};\xref,-\zr)}
 \Psi  (x^{\rm S},z^{\rm S};\xref,-\zr)\big) \Big]  ,
 \label{eq:KMmean1}
\end{eqnarray}
where, for any integer $j$, we have defined
\begin{eqnarray}
\label{def:Phij}
\Phi_{j}(x) &=& \frac{1}{N} \sum_{n=1}^N \beta_n^{j} \phi_n^2(x),\\
\label{def:Psi}
\Psi (x^{\rm S},z^{\rm S};\xref,\zr) &=&\frac{1}{N} \sum_{n=1}^N   \phi_n(\xref) \phi_n(x^{\rm S}) e^{i \beta_n (\zr-z^{\rm S})} .
\end{eqnarray}

The first term in the right-hand side of (\ref{eq:KMmean1})  is the contribution of the primary field.
The second term is the contribution of the secondary field.
The third term is a crossed contribution.

These results show that, when the waveguide is randomly perturbed and long enough (longer than the equipartition distance), then 
the illumination of the reflector becomes incoherent and Kirchhoff migration, which is based on coherent effects, 
gives a completely unstable and noisy image.

The analysis is complete in the time-harmonic case. 
The analysis of the broadband case (when the support of the source spectrum is not reduced to a single carrier frequency)
goes along the same line although it is necessary to use the asymptotic expressions 
of the two-frequency second-order moments of the transmission coefficients
(see \cite[Proposition 20.7]{FGPS} or \cite[Proposition 6.3]{GP-SIAM07}): due to the loss of coherence, the mean of the imaging functional is zero while its
variance is not.

\section{Correlation-Based Imaging Functionals}
\label{sec:imagf}

In this section, we introduce a new imaging functional to localize the point reflector using the scalar field recorded at the receiver array at $z=\widetilde{L}$. This functional is based on the correlation functions of the recorded signals, which we stress in the first subsection.

\subsection{Correlation of the scalar field}

Let $\sA$ denote the positions of the receivers in the plane $z=\widetilde{L}$. The data of scalar (pressure) field recorded by the receivers are $\{p^\eps(t,x,\widetilde{L}) ~|~ t\in \R, x \in \sA\}$. For simplicity, we have assumed that the receivers record data at all time. From these data one can form the 
cross correlation of the recorded field:
\begin{equation}
\Cc (\tau, x_1, x_2) = \int_\R \overline{p}(t,x_1,\widetilde{L}) 
p(t+\tau,x_2,\widetilde{L}) dt, \quad
x_1,x_2 \in \sA.
\label{eq:Ccdef}
\end{equation}
In Fourier domain, it has the form:
\begin{equation}
\Cc(\tau, x_1, x_2) = \frac{1}{2\pi} \int_\R \overline{\hat{p}}
(\omega, x_1, \widetilde{L}) \hat{p} (\omega, x_2, 
\widetilde{L}) e^{-i\omega \tau} d\omega, \quad
x_1,x_2 \in \sA.
\label{eq:Ccdef2}
\end{equation}

Using the decomposition $\hat{p} = \hat{p}_{\rm p}  + \hat{p}_{\rm s} $ in \eqref{eq:phatdec}, we can decompose the above cross correlation function into four parts. Let $\Cc_{\rm pp}$ denote the cross correlation between the primary fields at the two receivers. Thanks to the formula \eqref{eq:ppfs}, it admits the expression 
\begin{equation}
\begin{aligned}
\Cc_{\rm pp}(\tau,x_1,x_2) =& \frac{1}{8\pi} \int 
\sum_{j,l,m,n=1}^{N(\omega)} 
\frac{1}{\sqrt{\beta_l \beta_n \beta_j \beta_m(\omega)}} \overline{T^\eps_{jl}}(\omega) T^\eps_{mn}(\omega) 
|\hat{f}(\omega)|^2 \phi_j(x_1) \phi_m(x_2) \\
& \phi_l(\xs) \phi_n(\xs) e^{i(\beta_m - \beta_j) \widetilde{L}} 
e^{-i\omega \tau} d\omega.
\end{aligned}
\label{eq:Ccppdef}
\end{equation}

Let $\Cc_{\rm ps}$ denote the cross correlation between the primary field at the first receiver with the secondary field at the second receiver. Recall that the secondary field contains information about the waves emitted from the reflector at $\mathbf{x}_{\rm r} = (\xref,\widetilde{z}_{\rm r})$,
with $\widetilde{z}_{\rm r} = \widetilde{L} +\zr$.
Due to \eqref{eq:ppfs} and \eqref{eq:psfs}, it admits the expression 
\begin{equation}
\begin{aligned}
\Cc_{\rm ps}(\tau,x_1,x_2) =&  \int \sum_{q,j,l,m,n=1}^{N(\omega)} 
\frac{i\omega^2 \sigref}{16\pi\beta_q  \sqrt{\beta_l \beta_n \beta_j \beta_m(\omega)}} 
 \overline{T^\eps_{jl}}(\omega) T^\eps_{mn}(\omega) 
|\hat{f}(\omega)|^2 \phi_j(x_1) \phi_q(x_2) \\
& \phi_l(\xs) \phi_n(\xs) \phi_q(\xref) \phi_m(\xref) e^{i(\beta_m - \beta_j) 
\widetilde{L}} e^{i(\beta_q + \beta_m) \zr} 
e^{-i\omega \tau} d\omega.
\end{aligned}
\label{eq:Ccpsdef}
\end{equation}

Similarly, let $\Cc_{\rm sp}$ denote the cross correlation between the secondary field at the first receiver with the primary field at the second receiver. One verifies that
\begin{equation}
\begin{aligned}
\Cc_{\rm sp}(\tau,x_1,x_2) =& \int \sum_{q,j,l,m,n=1}^{N(\omega)} 
\frac{-i\omega^2 \sigref}{16\pi \beta_q \sqrt{\beta_l \beta_n \beta_j \beta_m(\omega)}} 
\overline{T^\eps_{jl}}(\omega) T^\eps_{mn}(\omega) 
|\hat{f}(\omega)|^2 \phi_q(x_1) \phi_m(x_2) \\
& \phi_l(\xs) \phi_n(\xs) \phi_q(\xref) \phi_j(\xref) e^{i(\beta_m - \beta_j) 
\widetilde{L}} e^{-i(\beta_q + \beta_j)\zr
} 
e^{-i\omega \tau} d\omega.
\end{aligned}
\label{eq:Ccspdef}
\end{equation}

Finally, the cross correlation between the secondary fields at the two receivers is much smaller than those above and its contribution is ignored. We neglected also the contributions from the error terms of the decomposition \eqref{eq:phatdec}. These are justified because $\sigref \ll 1$ consistently with the Born approximation.

Recall that the source in the acoustic model is due to the force $F(t,\mathbf{x})=f(t)\delta(\mathbf{x}-\mathbf{x}_{\mathrm{source}})$ where $\mathbf{x}_{\mathrm{source}}=(\xs,0)$ indicates the location of the source. In the rest of the paper, we will consider two special cases as follows. 

\subsubsection{Cross correlation for broadband pulse}
\label{subsec:crossbroad}%
We first consider the case where the source is given by $F(t,\mathbf{x})=f (t)\delta(\mathbf{x}-\mathbf{x}_{\mathrm{source}})$ with
\begin{equation}
\label{broadbandsource1}
f (t) = f_0(\eps^\alpha t) e^{-i\omega_0 t}.
\end{equation}
Here, $\omega_0$ is the carrier frequency. In the Fourier domain, we have
\begin{equation*}
\hat{f} (\omega) = \frac{1}{\eps^\alpha} \hat{f}_0
\Big( \frac{\omega-\omega_0}{\eps^\alpha} \Big).
\end{equation*}
Here, $\hat{f}_0$ is assumed to be a function with 
compact support or fast decay.

%[REVISION: ALPHA=2 BELONGS TO NARROWBAND: changed $\alpha>2$ to $\alpha\ge 2$]
When $\alpha \ge 2$ the bandwidth has no effect and the situation is equivalent to the time-harmonic case that we address in the next section.
%[END OF REVISION]

When $\alpha = (0,2)$, the pulse is said to be broadband and the bandwidth plays a role in the propagation
in the waveguide for a propagation distance of the order of $\widetilde{L}=L/\eps^2$.
Although the analysis can be carried out in general, we restrict ourselves
 to the case $\alpha \in (1,2)$ because when $\alpha \leq 1$, the number of propagating modes $N(\omega)$ varies with $\omega$ over the bandwidth
and the analysis is a little bit more delicate. 
Nevertheless, the overall picture does not change in the latter case.

Henceforth, $\alpha$ is a fixed number in the interval $(1,2)$. Let $\omega = \omega_0 + \eps^\alpha h$. Then
\begin{equation*}
\hat{f} (\omega) = \frac{1}{\eps^\alpha} \hat{f}_0(h), \quad T^\eps_{jl}(\omega) = 
T^\eps_{jl}(\omega_0 + \eps^\alpha h),
\end{equation*}
in terms of the new variable $h$. Further, we have the following Taylor expansions
\begin{equation*}
 \beta_j(\omega) =  \beta_j + \eps^\alpha \beta'_j h + 
o(\eps^\alpha), \quad
\frac{1}{\sqrt{\beta_j \beta_m \beta_l \beta_n(\omega)}} = \frac{1}{\sqrt{\beta_j \beta_m \beta_l \beta_n}} 
+ O(\eps^\alpha).
\end{equation*}
Here, $\beta'_j$ is the derivative of $\beta_j$ at the carrier frequency $\omega_0$; further, the reduced wavenumber $\beta_j$ is also evaluated at $\omega_0$. Using these formulas, the cross correlation functions become
\begin{equation*}
\begin{aligned}
\Cc_{\rm pp}(\tau,x_1,x_2) \approx & \frac{1}{8\pi\eps^\alpha} \int 
\sum_{j,l,m,n=1}^{N} \frac{|\hat{f}_0(h)|^2}{\sqrt{\beta_j \beta_m \beta_l \beta_n}} 
\overline{T^\eps_{jl}}(\omega_0+\eps^\alpha h) T^\eps_{mn}(\omega_0 + 
\eps^\alpha h)  \phi_j(x_1) \phi_m(x_2) \\
& \phi_l(\xs) \phi_n(\xs) e^{i[\beta_m (\omega_0 + \eps^\alpha h) 
-\beta_j(\omega_0 + \eps^\alpha h)]\widetilde{L}}  e^{-i(\omega_0 
+\eps^\alpha h)\tau} dh ,
\end{aligned}
\end{equation*}
\begin{equation*}
\begin{aligned}
\Cc_{\rm ps} \approx&  \int \sum_{q,j,l,m,n=1}^{N} \frac{i\omega_0^2 
\sigref}{16\pi\beta_q \eps^\alpha} 
\frac{|\hat{f}_0(h)|^2}{\sqrt{\beta_j \beta_m \beta_l \beta_n}} 
\overline{T^\eps_{jl}}(\omega_0 + \eps^\alpha h) T^\eps_{mn}(\omega_0 + 
\eps^\alpha h)  \phi_j(x_1) \phi_q(x_2) \phi_l(\xs) \\
&  \phi_n(\xs)\phi_q(\xref) \phi_m(\xref) e^{i[\beta_m (\omega_0 + 
\eps^\alpha h) -\beta_j(\omega_0 + \eps^\alpha h)]\widetilde{L}} 
e^{i(\beta_q + \beta_m)\zr} e^{i(\beta'_q + \beta'_m)\eps^\alpha h \zr} 
e^{-i(\omega_0+\eps^\alpha h) 
\tau} dh ,
\end{aligned}
\end{equation*}
\begin{equation*}
\begin{aligned}
\Cc_{\rm sp} \approx & \int \sum_{q,j,l,m,n=1}^{N} \frac{-i\omega_0^2 
\sigref}{16\pi \beta_q\eps^\alpha} 
\frac{|\hat{f}_0(h)|^2}{\sqrt{\beta_j \beta_m \beta_l \beta_n}} 
\overline{T^\eps_{jl}}(\omega_0 + \eps^\alpha h) T^\eps_{mn}(\omega_0 + 
\eps^\alpha h)  \phi_q(x_1) \phi_m(x_2) \phi_l(\xs) \\
& \phi_n(\xs) \phi_q(\xref) \phi_j(\xref) e^{i[\beta_m (\omega_0 + 
\eps^\alpha h) -\beta_j(\omega_0 + \eps^\alpha h)] 
\widetilde{L}}e^{-i(\beta_q + \beta_j)\zr} 
e^{-i(\beta'_q + \beta'_j) \eps^\alpha h \zr } 
e^{-i(\omega_0+\eps^\alpha h) \tau} dh.
\end{aligned}
\end{equation*}

%%%%%%%%%%%%%%%%%%
\subsubsection{Cross correlation with time-harmonic source}

A simple case is when the source term is time-harmonic, {\itshape i.e.} $F(t,\mathbf{x})=\delta(\mathbf{x}-\mathbf{x}_{\mathrm{source}}) f(t) $
with $f(t) = e^{-i\omega_0 t}$ and
\begin{equation*}
\hat{f}(\omega) = 2\pi \delta(\omega-\omega_0).
\end{equation*}
%[REVISION]
%[REMOVED]Formally, one substitutes the above relation into the expressions \eqref{eq:Ccppdef}, \eqref{eq:Ccpsdef} and \eqref{eq:Ccspdef} and observes that the cross correlations are simplified to[REMOVED]
In this case, the wave field has the form $p(t,x,z) = \hat{p}(x,z;\omega_0) e^{-i\omega_0 t}$. The definition of the correlation function should be modified to
\begin{equation}
\Cc(\tau,x_1,x_2) : = \frac{1}{T} \int_0^T 
\overline{p}(t,x_1,\widetilde{L}) p(t+\tau, x_2, \widetilde{L}) dt 
= e^{-i\omega_0\tau} \overline{\hat{p}}(x_1,\widetilde{L};\omega_0) 
\hat{p}(x_2,\widetilde{L};\omega_0).
\label{eq:Ccsf}
\end{equation}
The second equality holds because the integrand above is in fact independent of $t$. Using the 
decomposition \eqref{eq:phatdec} and the expressions \eqref{eq:ppfs} and \eqref{eq:psfs}, we obtain the following expressions for the cross correlations:
%[END OF REVISION]
\begin{equation*}
\begin{aligned}
\Cc_{\rm pp}(\tau,x_1,x_2) =& \frac{1}{4} \sum_{j,l,m,n=1}^{N(\omega_0)} 
\frac{1}{\sqrt{\beta_j \beta_m \beta_l\beta_n}} 
\overline{T^\eps_{jl}}(\omega_0) T^\eps_{mn}(\omega_0) \phi_j(x_1) 
\phi_m(x_2) \\
& \phi_l(\xs) \phi_n(\xs) e^{i(\beta_m - \beta_j) \widetilde{L}} 
e^{-i\omega_0 \tau},
\end{aligned}
\end{equation*}

\begin{equation*}
\begin{aligned}
\Cc_{\rm ps}(\tau,x_1,x_2) =&  \sum_{q,j,l,m,n=1}^{N(\omega_0)} 
\frac{i\omega_0^2 \sigref}{8\beta_q} 
\frac{1}{\sqrt{\beta_j \beta_m \beta_l\beta_n}} 
\overline{T^\eps_{jl}}(\omega_0) T^\eps_{mn}(\omega_0) \phi_j(x_1) 
\phi_q(x_2) \\
& \phi_l(\xs) \phi_n(\xs) \phi_q(\xref) \phi_m(\xref) e^{i(\beta_m - \beta_j) 
\widetilde{L}} e^{i(\beta_q + \beta_m)\zr
} 
e^{-i\omega_0 \tau},
\end{aligned}
\end{equation*}

\begin{equation*}
\begin{aligned}
\Cc_{\rm sp}(\tau,x_1,x_2) =& \sum_{q,j,l,m,n=1}^{N(\omega_0)} 
\frac{-i\omega_0^2 \sigref}{8 \beta_q} 
\frac{1}{\sqrt{\beta_j \beta_m \beta_l\beta_n}} 
\overline{T^\eps_{jl}}(\omega_0) T^\eps_{mn}(\omega_0) \phi_q(x_1) 
\phi_m(x_2) \\
& \phi_l(\xs) \phi_n(\xs) \phi_q(\xref) \phi_j(\xref) e^{i(\beta_m - \beta_j) 
\widetilde{L}} e^{-i(\beta_q + \beta_j)\zr} 
e^{-i\omega_0 \tau}.
\end{aligned}
\end{equation*}

%%%%%%%%%%%%%%%
\subsection{Imaging functionals using cross correlations}

We are now ready to present the imaging functionals,
which consist in migrating the cross correlations of the recorded signals.
The imaging functionals are designed according to the settings of receiver arrays. We consider two cases.

\bigskip
{\bfseries Full aperture receiver array.} The ideal case is when the receiver array spans the whole
cross section, {\itshape i.e.}, $\sA = \mathcal{D}$. Then given the data, for any pair of modes $\phi_j$ and $\phi_l$, we define
\begin{equation}
\Cc^{jl}(\tau) : = \int_{\sA} \int_{\sA} \Cc(\tau,x_1,x_2) \phi_j(x_1) \phi_l(x_2) dx_1 
dx_2.
\label{eq:Ccjldef}
\end{equation}
Due to the orthogonality of $\{\phi_j\}_{j=1,\ldots,N}$, the function $\Cc^{jl}$ is the $jl$ mode of the cross correlation function $\Cc$.

A search point for the reflector will be denoted as  $(x^{\rm S},\widetilde{z}^{\rm S})$ where $x^{\rm S}$ is its transversal coordinate
and $\widetilde{z}^{\rm S}=\widetilde{L} + z^{\rm S}$ is its axial coordinate. 
Equivalently, $z^{\rm S}$ is the axial coordinate starting from the receiver array. We define the imaging
 functional $\mathcal{I}_{\mathrm{FA}}$ as
\begin{eqnarray}
\mathcal{I}_{\rm FA} (x^{\rm S}, z^{\rm S}) &=&
\mathcal{I}_{{\rm FA}+} (x^{\rm S}, z^{\rm S}) +
\mathcal{I}_{{\rm FA}-} (x^{\rm S}, z^{\rm S}) ,\\
\mathcal{I}_{{\rm FA}\pm} (x^{\rm S}, z^{\rm S})&=& \mp 
\frac{i}{N(\omega_0)^2}\sum_{j,l=1}^{N(\omega_0)} \beta_j \beta_l (\omega_0)\phi_j(x^{\rm S}) \phi_l(x^{\rm S})  \Cc^{jl}
\Big(\pm \frac{z^{\rm S}}{\omega_0} (\beta_j + \beta_l)\Big) .
\label{eq:IKMdef}
\end{eqnarray}
The choice of multiplication by $\beta_j \beta_l$ is suggested from our analysis of the correlation functions in the next section.

\bigskip
{\bfseries Limited aperture receiver array.} A realistic situation is when the receiver array only covers part of the transversal section, {\itshape i.e.}, $\sA = [a_1,a_2]$ for $0< a_1 < a_2 <a$. Consequently, the exact $jl$ component of the cross correlation cannot be extracted. In this case, we design the following imaging functional $\mathcal{I}_{\rm LA}$ as
\begin{eqnarray}
\mathcal{I}_{\rm LA} (x^{\rm S}, z^{\rm S}) &=&
\mathcal{I}_{{\rm LA}+} (x^{\rm S}, z^{\rm S}) +
\mathcal{I}_{{\rm LA}-} (x^{\rm S}, z^{\rm S}) ,\\
\mathcal{I}_{{\rm LA}\pm} (x^{\rm S},z^{\rm S}) &=&  \mp \frac{ i }{N(\omega_0)^2}\int_{\mathcal{A}^2} \sum_{j,q=1}^{N(\omega_0)} \phi_{q}(x_1) 
\phi_{q}(x^{\rm S}) \phi_{j}(x_2) \phi_{j}(x^{\rm S})
\label{eq:Idef}\\
&& \times (\Delta_{x_1} +k(\omega_0)^2)(\Delta_{x_2} + k(\omega_0)^2) 
\Cc \Big( \pm \frac{z^{\rm S}}{\omega_0}(\beta_{q} + 
\beta_{j}),x_1,x_2 \Big)  dx_1 dx_2.\notag
\end{eqnarray}
Again, these definitions of imaging functionals are based on the analysis of the correlation functions in the next section. We remark also that it is possible to show that when $\mathcal{A} = \mathcal{D}$, the second functional $\mathcal{I}_{{\rm LA}}$ is very close to $\mathcal{I}_{{\rm FA}}$ and we have in fact:
$$
\mathcal{I}_{{\rm LA}\pm} (x^{\rm S}, z^{\rm S}) \mid_{{\mathcal A}={\mathcal D}}= \mp 
\frac{i}{N(\omega_0)^2}\sum_{j,l=1}^{N(\omega_0)} \beta_j^2 \beta_l^2 (\omega_0)\phi_j(x^{\rm S}) \phi_l(x^{\rm S})  \Cc^{jl}
\Big(\pm \frac{z^{\rm S}}{\omega_0} (\beta_j + \beta_l)\Big) .
$$

%%%%%%%%%%%%%%%
%%%%%%%%%%%%%%%
\section{Resolution Analysis of the Imaging Functionals}
\label{sec:resol}

In this section, we analyze the imaging functionals proposed above to search the reflectors in the waveguide. Due to the random perturbations of the
long section $z \in [0,\widetilde{L}]$ of the waveguide, the values of the imaging functionals, which depend on the waveguide parameters through the data, are random. Hence, we analyze the mean of the imaging functional and show that it achieves its maximum at the reflector location. We also analyze how this mean decays from its maximum; this information provides the resolution of the proposed imaging functionals.

We emphasize that the above observation on the mean of the imaging functional itself is not enough to claim that the functionals are effective, because it is not certain, {\itshape a priori}, that the one realization in practice is well reflected by the mean. Statistical stability
({\itshape i.e.}, the analysis of the fluctuations of the imaging functionals) is needed to secure this claim. This will be investigated in the Section \ref{sec:stabi}.

%%%%%%%%%%%%%%%
\subsection{The case of full aperture receiver array with time-harmonic sources}

We first consider the ideal case where the receiver array $\sA$ spans the whole cross section, 
so the imaging functional $\mathcal{I}_{\mathrm{FA}}$ is chosen. 
The key tool for analysis of the mean of $\mathcal{I}_{\mathrm{FA}}$ is Proposition \ref{prop:momentT}.

As shown by Proposition \ref{prop:momentT}, when $L$ is larger than the energy equipartition length 
$L_{\mathrm{equip}}$, the main contribution from terms of 
$\E\{T^\eps_{jl}T^\eps_{mn}\}$ comes from when $j=m$ and $l=n$.  
Following this observation, we have the following limits for the cross 
correlations.

\begin{equation*}
\begin{aligned}
\E\{\Cc_{\rm pp}(\tau,x_1,x_2)\}  
& \rightarrow \frac{\Phi_{-1}(x_{\rm s})}{4} \sum_{j=1}^N \frac{1}{\beta_j} \phi_j(x_1) 
\phi_j(x_2) e^{-i\omega_0 \tau} ,
\end{aligned}
\end{equation*}
where $\Phi_{-1}$ is defined by (\ref{def:Phij}) and we have also used the fact that $\mathcal{T}^{(l)}_j$ converges to 
$1/N$ in the regime $L \gg L_{\mathrm{equip}}$. 

Following the same lines, we have
\begin{equation}
\E\{\Cc_{\rm ps}(\tau,x_1,x_2)\} \rightarrow \frac{i\omega_0^2 \sigref \Phi_{-1}(x_{\rm s})}{8} 
\sum_{q,j=1}^{N(\omega_0)} \frac{1}{\beta_j \beta_q} \phi_j(x_1) 
\phi_j(\xref) \phi_q(\xref) \phi_q(x_2) 
e^{i(\beta_j+\beta_q)\zr} e^{-i\omega_0 \tau} ,
\label{eq:Ccpssf}
\end{equation}
and
\begin{equation}
\E\{\Cc_{\rm sp}(\tau,x_1,x_2)\} \rightarrow \frac{-i\omega_0^2 \sigref \Phi_{-1}(x_{\rm s})}{8} \sum_{q,j=1}^{N(\omega_0)} \frac{1}{\beta_j \beta_q} \phi_q(x_1) 
\phi_q(\xref) \phi_j(\xref) \phi_j(x_2) 
e^{-i(\beta_q+\beta_j)\zr} e^{-i\omega_0 \tau}.
\label{eq:Ccspsf}
\end{equation}
In fact, the migration imaging functional ${\mathcal I}_{{\rm FA}}$ is designed from the above 
characterization of the cross correlation function. From the calculations before, we find
\begin{equation*}
\E \{\Cc^{jl}_{\rm pp}(\tau)\} \rightarrow \frac{\Phi_{-1}(x_{\rm s})}{4\beta_j} \delta_{jl} 
e^{-i\omega_0 \tau},
\end{equation*}
where $\delta_{jl}$ is the Kronecker symbol. We also find
\begin{equation*}
\E \{\Cc^{jl}_{\rm ps}(\tau)\} \rightarrow \frac{i\omega_0^2 \sigref \Phi_{-1}(x_{\rm s})}{8} 
\frac{1}{\beta_j \beta_l}  \phi_j(\xref) \phi_l(\xref) 
e^{i(\beta_j+\beta_l)\zr} e^{-i\omega_0 \tau}.
\end{equation*}
\begin{equation*}
\E \{\Cc^{jl}_{\rm sp}(\tau)\} \rightarrow \frac{-i\omega_0^2 \sigref \Phi_{-1}(x_{\rm s})}{8} 
\frac{1}{\beta_j \beta_l}  \phi_j(\xref) \phi_l(\xref) 
e^{-i(\beta_j+\beta_l)\zr} e^{-i\omega_0 \tau}.
\end{equation*}

Therefore, for a search point $(x^{\rm S},\widetilde{z}^{\rm S})$, 
with $\widetilde{z}^{\rm S}=\widetilde{L}+z^{\rm S}$, we have the following. From now on, $\lambda=2\pi/k(\omega_0)$ denotes the 
carrier wavelength.
\begin{proposition}
If $a \gg \lambda$, $z^{\rm S},\zr \gg \lambda$, and $x^{\rm S},\xref \in (0,a)$, then
\begin{equation}
\E \big[ \mathcal{I}_{{\rm FA}} (x^{\rm S}, z^{\rm S}) \big] \simeq \frac{\pi \omega_0^2 \sigref}{32a^3}
\Re e \Big\{ \Big[ \int_{-\frac{\pi}{2}}^{\frac{\pi}{2}} \cos\theta 
e^{i(\tilde{\eta} \cos\theta + \tilde{\xi} \sin\theta)}d\theta
\Big]^2 \Big\} ,
\label{eq:prop1}
\end{equation}
where we have introduced the normalized cross 
range offset $\tilde{\xi} = 2\pi(\xref - x^{\rm S}
)/\lambda$, and the normalized 
range offset $\tilde{\eta} = 2\pi(\zr - z^{\rm S}
)/\lambda$. 
\end{proposition}

Therefore, the imaging functional $\mathcal{I}_{\mathrm{FA}}$ 
works well to detect the point reflector $(\xref,\widetilde{z}_{\rm r})$. In particular we can see that both the range 
and cross-range resolutions are of the order of  the wavelength $\lambda$.

\begin{proof}
Note that $a \gg \lambda$ means that $N \gg 1$.
Calculations show that
\begin{eqnarray}
\nonumber
\E \big[ \mathcal{I}_{{\rm FA}+} (x^{\rm S},z^{\rm S}) \big] 
&\rightarrow & \frac{-i \Phi_{-1}(x_{\rm s})}{4N}
\Big[ \frac{1}{N}  \sum_{j=1}^N \beta_j \phi^2_j(x^{\rm S}) e^{-i2\beta_j z^{\rm S}} \Big] \\
&& \hspace*{-0.5in} +  \frac{\omega_0^2 \sigref \Phi_{-1}(x_{\rm s})}{8} 
\Psi(x^{\rm S},z^{\rm S};x_{\rm r},z_{\rm r})^2
- 
 \frac{\omega_0^2 \sigref \Phi_{-1}(x_{\rm s})}{8}
\Psi(x^{\rm S},z^{\rm S};x_{\rm r},-z_{\rm r})^2  ,
\hspace*{0.2in}
\label{eq:espIFAp}
\end{eqnarray}
where $\Psi$ and $\Phi_j$ are defined by (\ref{def:Psi}) and (\ref{def:Phij}).
The first term does not contain information about the reflector and can 
be viewed as background field. In fact, its contribution is negligible because the fast oscillations in the complex exponential ($\zr , z^{\rm S} \gg \lambda$) average out and there is an overall factor $1/N$. For the second term, we can use an integral 
approximation for the sum in the continuum limit ($N \gg 1$). That is,
\begin{equation*}
\begin{aligned}
& 
\Psi(x^{\rm S},z^{\rm S};x_{\rm r},z_{\rm r})
=  \frac{1}{N} \sum_{j=1}^N \frac{2}{a} \sin(\frac{2\pi x^{\rm S}}{\lambda} 
\frac{j}{N}) \sin(\frac{2\pi \xref}{\lambda} \frac{j}{N}) e^{i2\pi 
\sqrt{1-(\frac{j}{N})^2} \frac{(\zr - z^{\rm S}
)}{\lambda}}\\
& = \frac{1}{aN} \sum_{j=1}^N \left(\cos(\frac{2\pi (\xref-x^{\rm S}
)}{\lambda} 
\frac{j}{N}) - \cos(\frac{2\pi (\xref+ x^{\rm S}
)}{\lambda} \frac{j}{N})\right)  
e^{i2\pi \sqrt{1-(\frac{j}{N})^2} \frac{\zr - z^{\rm S}
}{\lambda}} \\
&\approx  \frac{1}{a} \int_0^1 \left( \cos(\frac{2\pi(\xref - 
x^{\rm S}
)}{\lambda} y) - \cos(\frac{2\pi(\xref + x^{\rm S}
)}{\lambda} y)\right) 
e^{i2\pi \sqrt{1-y^2} \frac{\zr - z^{\rm S}
}{\lambda}} dy.
\end{aligned}
\end{equation*}
Since the phase becomes zero when $z^{\rm S}
 = \zr$, this function peaks at $z^{\rm S}= \zr$ and $x^{\rm S}= \xref$.
The integral 
considered above can be written as, with the second term neglected,
\begin{equation*}
\begin{aligned}
\frac{1}{2a} \int_0^1 e^{i(\tilde{\eta} \sqrt{1-y^2} + \tilde{\xi} y)} + 
e^{i(\tilde{\eta} \sqrt{1-y^2} - \tilde{\xi} y)} dy  &= \frac{1}{2a} 
\int_{-1}^1 e^{i(\tilde{\eta} \sqrt{1-y^2} + \tilde{\xi} y)} dy\\
&= \frac{1}{2a} \int_{-\frac{\pi}{2}}^{\frac{\pi}{2}} \cos\theta 
e^{i(\tilde{\eta} \cos\theta + \tilde{\xi} \sin\theta)}d\theta.
\end{aligned}
\end{equation*}

Similarly, the sum in the third term can be approximated by
\begin{equation*}
\begin{aligned}
&
\Psi(x^{\rm S},z^{\rm S};x_{\rm r}, -z_{\rm r})
\approx  \frac{1}{a} \int_0^1 \left( \cos(\frac{2\pi(\xref - 
x^{\rm S}
)}{\lambda} y) - \cos(\frac{2\pi(\xref + x^{\rm S}
)}{\lambda} y)\right) 
e^{-i2\pi \sqrt{1-y^2} \frac{\zr + z^{\rm S}
}{\lambda}} dy.
\end{aligned}
\end{equation*}
Note that this function does not have a peak comparable with the previous 
function. In fact, stationary phase calculation shows that it is of order $O(1/(a\sqrt{k(\zr + z^{\rm S})}))$ where $k = 2\pi/\lambda$.

The  evaluation of  $\E [ \mathcal{I}_{\mathrm{FA}-}]$ follows the same lines:
\begin{eqnarray}
\nonumber
\E \big[ \mathcal{I}_{{\rm FA}-} (x^{\rm S},z^{\rm S}) \big] 
&\rightarrow & \frac{i \Phi_{-1}(x_{\rm s})}{4N}
\Big[ \frac{1}{N}  \sum_{j=1}^N \beta_j \phi^2_j(x^{\rm S}) e^{i2\beta_j z^{\rm S}} \Big] \\
&& \hspace*{-0.5in} -  \frac{\omega_0^2 \sigref \Phi_{-1}(x_{\rm s})}{8} 
\Psi(x^{\rm S},-z^{\rm S};x_{\rm r},z_{\rm r})^2
+
 \frac{\omega_0^2 \sigref \Phi_{-1}(x_{\rm s})}{8}
\Psi(x^{\rm S},-z^{\rm S};x_{\rm r},-z_{\rm r})^2  ,
\hspace*{0.2in}
\label{eq:espIFAm}
\end{eqnarray}
 and we find a peak at $z^{\rm S}
 = \zr$ and $x^{\rm S}= \xref$. 
 
 Finally we have
\begin{equation}
\label{asy:Phi-1}
 \Phi_{-1}(x_{\rm s}) = \frac{2}{a N} \sum_{j=1}^{N} \beta_j^{-1}  
  \sin^2(\frac{2\pi x_{\rm s}}{\lambda} 
\frac{j}{N})  \stackrel{N \gg 1}{\simeq} \frac{1}{a} \int_0^1 \frac{1}{\sqrt{1-x^2}} dx = \frac{\pi}{2a},
\end{equation}
which completes the proof of the proposition.
\end{proof}

\begin{remark} 
To eliminate the $\beta_j$ on the denominator, we multiplied by the exact 
$\beta_j$ factor in constructing the imaging function. Alternatively, we 
can differentiate the cross correlation function to remove the 
denominator; see \eqref{eq:Idef}. If we do so, the above approximation will lead to 
the integral
\begin{equation*}
 \int_{-\frac{\pi}{2}}^{\frac{\pi}{2}} \cos^2 \theta 
e^{i(\tilde{\eta} \cos\theta + \tilde{\xi} \sin\theta)}d\theta
\end{equation*}
in the square brackets of Eq.~(\ref{eq:prop1}).
This function should be compared with \eqref{eq:EILA} in the next section.  
\end{remark}

%%%%%%%%%%%%%%%
\subsection{The case of limited aperture receiver array}

Next, we consider the realistic setting where the receiver array only covers part of the transverse section, namely $\sA = (a_1, a_2)$ and $0<a_1 < a_2 <a$. For a search point $(x^{\rm S},\widetilde{z}^{\rm S})$, 
with $\widetilde{z}^{\rm S} = \widetilde{L} + {z}^{\rm S}$, we have the following result.
\begin{proposition}
\label{prop2}%
If $a_2-a_1 \gg \lambda$, $z^{\rm S},\zr \gg a$, and $x^{\rm S},\xref \in (0,a)$, then
\begin{equation}
\E \big[ \mathcal{I}_{{\rm LA}}  (x^{\rm S}, z^{\rm S}) \big] \simeq \frac{\pi \omega_0^2 k^2 \sigref (a_2-a_1)^2}{32a^5}
\Re e \Big\{ \Big[ \int_{-\frac{\pi}{2}}^{\frac{\pi}{2}} \cos^2\theta 
e^{i(\tilde{\eta} \cos\theta + \tilde{\xi} \sin\theta)}d\theta
\Big]^2 \Big\},
\label{eq:EILA}
\end{equation}
where we have introduced the normalized cross 
range offset $\tilde{\xi} = 2\pi(\xref - x^{\rm S})/\lambda$, and the normalized 
range offset $\tilde{\eta} = 2\pi(\zr - z^{\rm S})/\lambda$. 
\end{proposition}

\begin{figure}
\begin{center}
\includegraphics[width=7cm]{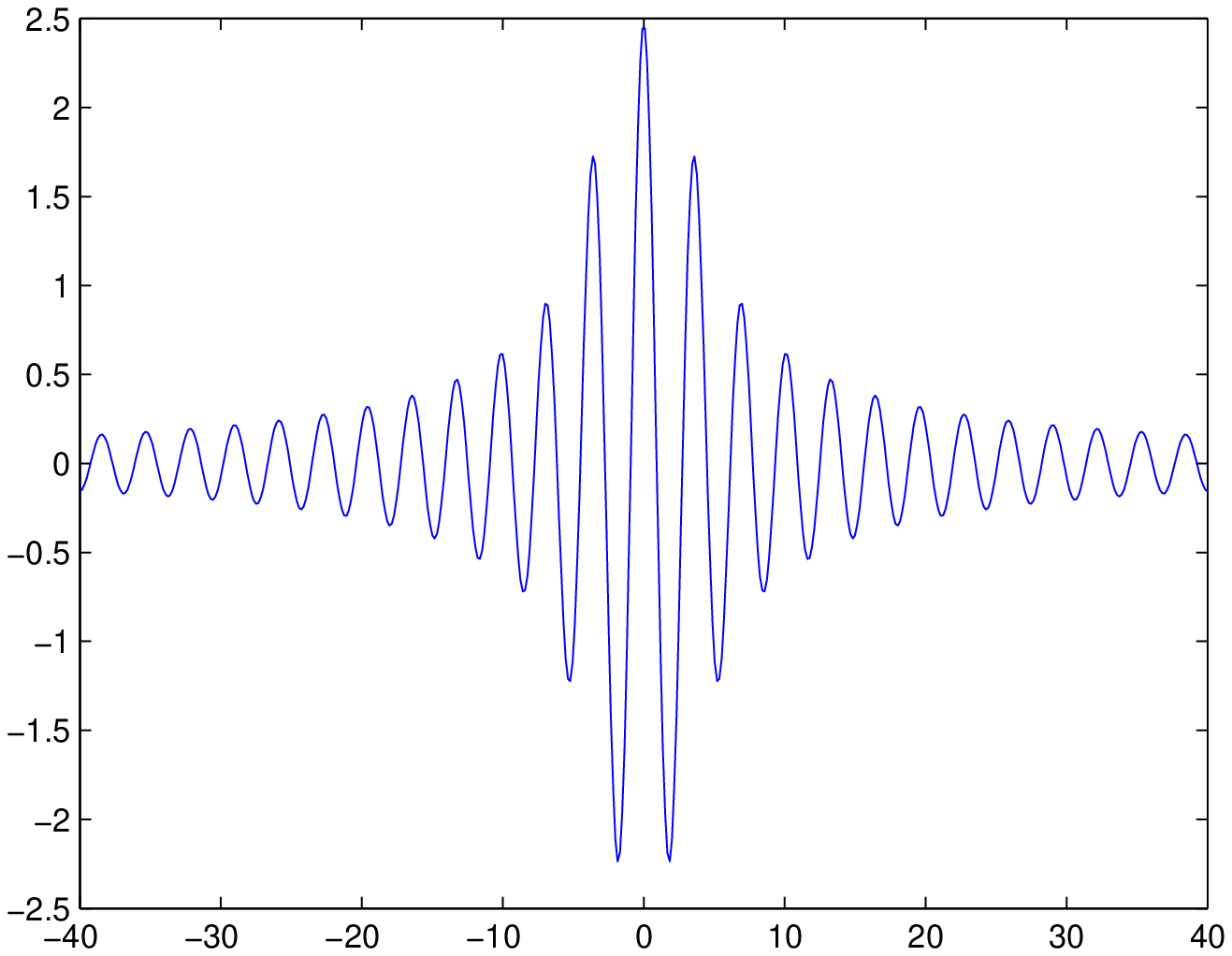}
\includegraphics[width=7cm]{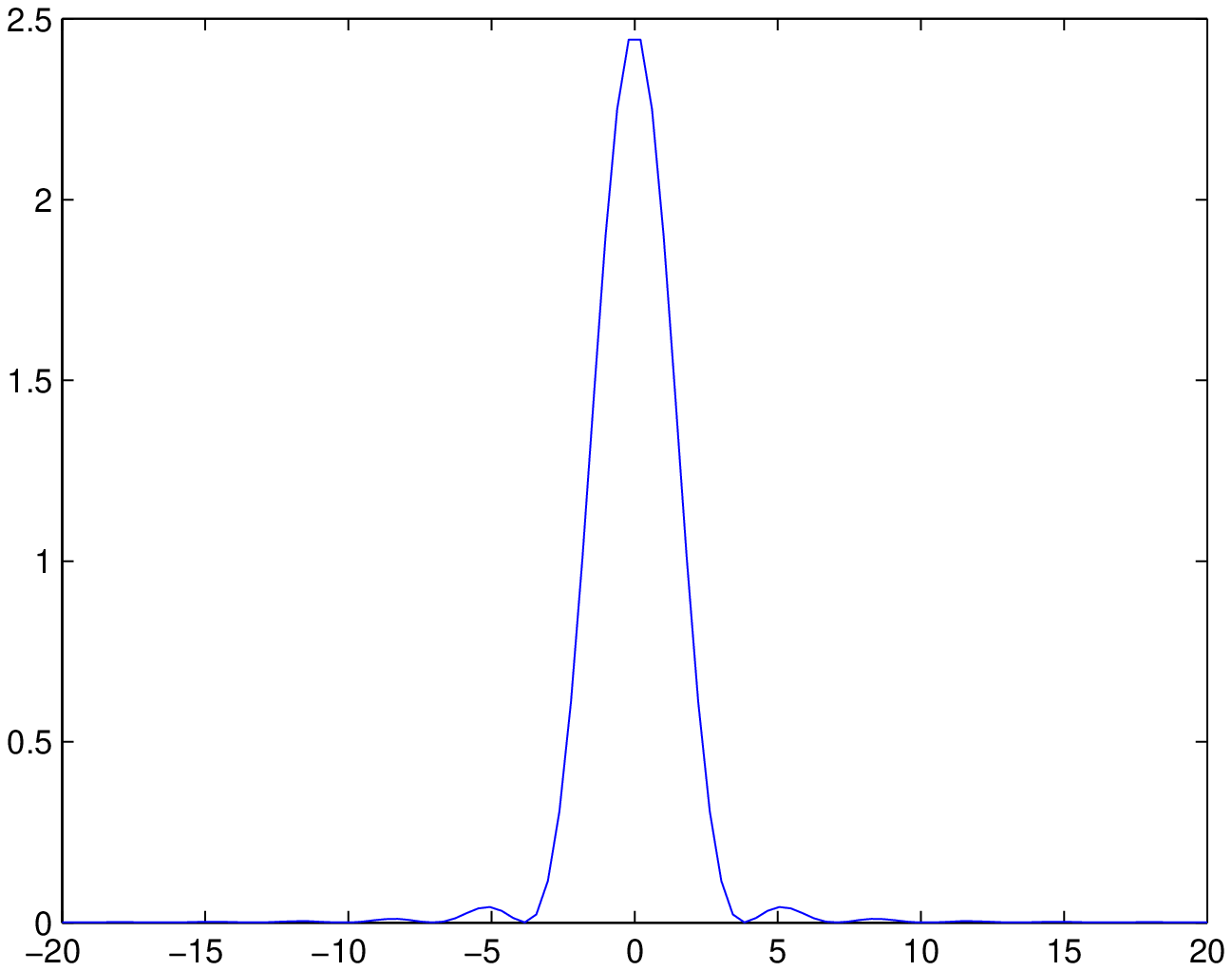}
\end{center}
\caption{Plots of the functions $g(\tilde{\eta})$ and $h(\tilde{\xi})$
which give the normalized form of the point spread function in the range direction ($g(\tilde{\eta})$, left picture)
and in the cross range direction ($h(\tilde{\xi})$, right picture).}
\label{fig1}
\end{figure}

This proposition gives the form of the point spread function of the mean imaging functional,
that is the normalized form of the peak centered at the reflector location.
The range and cross-range widths of the peak are the range and cross-range resolutions.
Since the variables $\tilde{\xi}$ and $\tilde{\eta}$ are normalized with 
respect to the wavelength, we see that both the range and cross-range 
resolutions are of order of the wavelength, which is the diffraction limit.
The form of the peak can can be seen in Figure \ref{fig1} which plots  the transverse and longitudinal shapes of the 
point spread function
\begin{eqnarray}
h(\tilde{\xi}) &=& \Re e \Big\{ \Big[ \int_{-\frac{\pi}{2}}^{\frac{\pi}{2}} \cos^2\theta e^{i  \tilde{\xi} \sin\theta}d\theta
\Big]^2 \Big\} = \pi^2 \frac{J_1(\tilde{\xi})^2}{\tilde{\xi}^2} , \\
g(\tilde{\eta})&=&  \Re e \Big\{ \Big[ \int_{-\frac{\pi}{2}}^{\frac{\pi}{2}} \cos^2\theta 
e^{i\tilde{\eta} \cos\theta }d\theta\Big]^2 \Big\} = \pi^2 J_1'( \tilde{\eta})^2 -   \Big[ \int_{-\frac{\pi}{2}}^{\frac{\pi}{2}} \cos^2\theta 
\sin( \tilde{\eta} \cos\theta )d\theta\Big]^2 ,
\end{eqnarray}
where we have used Formula 9.1.20 \cite{abra} to express $h$ and $g$ in terms of the 
Bessel function~$J_1$.
$h$ and $g$ are even functions maximal at $0$ and 
stationary phase calculations also show that
$$
h(\tilde{\xi}) \stackrel{\tilde{\xi} \gg 1}{\simeq}  \pi   \frac{1-\sin(2 \tilde{\xi})}{\tilde{\xi}^3} , \quad \quad
g(\tilde{\eta})\stackrel{\tilde{\eta} \gg 1}{\simeq }  2\pi\frac{ \sin(2\tilde{\eta} )}{\tilde{\eta}}  .
$$
Proposition \ref{prop2} also shows that the resolution of the mean imaging functional 
does not depend on the aperture of the array $a_2-a_1$. This is a consequence of the waveguide geometry, since multiple reflections at the boundaries of the waveguide generate multiple replicas of the receiver array in the plane $z=\widetilde{L}$, 
which gives an effective aperture that is  large enough to reach the diffraction limit.

\begin{proof}
Since the two functions above can be 
analyzed in the same manner, we focus on one of them, namely 
$\mathcal{I}_{{\rm LA}+}$.
As before, the cross correlation $\Cc$ can be treated component-wise. For the function $\mathcal{I}_{{\rm LA}+}$, the main contribution comes from the primary-secondary component $\Cc_{\rm ps}$, which we focus on for the moment. Let $\mathcal{I}_{+{\rm ps}}$ denote this main term, {\itshape i.e.}
\begin{eqnarray*}
\mathcal{I}_{+{\rm ps}}&=& - \frac{ i }{N^2}
\int_{\mathcal{A}^2} \sum_{j',q'=1}^N  \phi_{q'}(x_1) \phi_{q'}(\xref^s)  \phi_{j'}(x_2) 
\phi_{j'}(\xref^s) \\
&&\times (\Delta_{x_1} + k^2) (\Delta_{x_2} + k^2) \Cc_{\rm ps} 
\Big(\frac{z^{\rm S}}{\omega_0}(\beta_{q'} + \beta_{j'}), x_1,x_2\Big)  dx_1 dx_2.
\end{eqnarray*}
Then from \eqref{eq:Ccpssf} and the fact that $\Delta_{x_1} \phi_j = -\lambda_j^2 \phi_j$ and $k^2 = \lambda_j^2 + \beta_j^2$, we verify that
\begin{equation*}
\begin{aligned}
\E [\mathcal{I}_{+{\rm ps}}] \longrightarrow &\frac{\omega_0^2 \sigref \Phi_{-1}(x_{\rm s})}{8N^2} \sum_{q,j,q',j'=1}^N 
\beta_j \beta_q \phi_j(\xref) \phi_{j'}(x^{\rm S}
) 
\phi_q(\xref) \phi_{q'}(x^{\rm S}
) e^{i(\beta_j + \beta_q)\zr} 
e^{-i(\beta_{j'}+\beta_{q'}) z^{\rm S}
}\\
&\times \int_{\mathcal{A}^2} \phi_j(x_1) \phi_{j'}(x_1) \phi_q(x_2) 
\phi_{q'}(x_2) dx_1 dx_2.
\end{aligned}
\end{equation*}
If $\mathcal{A} = (0,a)$, the last integral will be $\delta_{jj'} 
\delta_{qq'}$ and we are back in the case of full aperture array. In the 
current situation, this integral has to be dealt with more carefully.  
The above limit can be further written as
\begin{equation}
\frac{\omega_0^2 \sigref \Phi_{-1}(x_{\rm s})}{8} \Big[ \frac{1}{N}\sum_{j,j'=1}^N 
\beta_j \phi_j(\xref) \phi_{j'}(x^{\rm S}
) e^{i(\beta_j \zr - 
\beta_{j'} z^{\rm S}
)} \int_{\mathcal{A}} \phi_j(x_1) \phi_{j'}(x_1)dx_1 
\Big]^2.
\label{eq:si1ps}
\end{equation}
Introducing the difference index $l = j - j'$, the double sum above can written in terms of $j'$ and $l$. The integral over $x_1$ can be calculated explicitly as follows:
\begin{equation*}
\begin{aligned}
\int_{\mathcal{A}} \phi_{j'+l}(x_1) \phi_{j'}(x_1) dx_1 = &\int_{a_1}^{a_2} 
\frac{2}{a} \sin(\frac{(j'+l) \pi x_1}{a}) \sin(\frac{j' \pi x_1}{a}) dx_1\\
= &\frac{1}{a} \int_{a_1}^{a_2} \cos( \frac{\pi l x_1}{a}) -  \cos( 
\frac{\pi(2j' + l) x_1}{a}) dx_1\\
=& \frac{a_2 - a_1}{a} \left[\cos(\frac{\pi l(a_2+a_1)}{2a}) \sinc
(\frac{\pi l (a_2 -a_1)}{2a})\right.\\
&\quad  - \left.\cos(\frac{\pi(2j' + l)(a_2+a_1)}{2a}) \sinc 
(\frac{\pi(2j' + l)(a_2 -a_1)}{2a})\right].
\end{aligned}
\end{equation*}
For  most $j'$, the first term above dominates. Hence, the number of indices $l$'s so that the above quantity is significant is roughly of order $a/(a_2 - a_1)$.

Using the explicit expression $\phi_j(x) = \sqrt{2/a} \sin(j\pi x/a)$, we rewrite the function inside the 
brackets in \eqref{eq:si1ps} as
\begin{equation*}
\begin{aligned}
 \int_{\mathcal{A}} \frac{1}{N}\sum_{j',l}
 \frac{4\beta_{j'+l}}{a^2} &\sin (\frac{\pi (j'+l) \xref}{a}) \sin (\frac{\pi j' x^{\rm S}
}{a})  \sin (\frac{\pi (j'+l) x_1}{a}) \sin (\frac{\pi j' x_1}{a})\\
& e^{i\beta_{j'} (\zr - z^{\rm S}
)} e^{i(\beta_{j'+l} - \beta_{j'}) \zr} dx_1.
\end{aligned}
\end{equation*}
Using the explicit expression $\beta_j = \sqrt{k^2 - (j\pi/a)^2}$ and the fact that the total number of modes $N$ is 
$\lfloor a/(\lambda/2) \rfloor = \lfloor ak/\pi \rfloor$ and assuming that $ak/\pi$ is an integer for simplicity, we see that for $l \ll N$,
\begin{equation}
\beta_{j'+l} - \beta_{j'} \approx -\frac{k^2}{\beta_{j'}} \frac{j'}{N} \frac{l}{N}.
\label{eq:betalin}
\end{equation}
Using this expansion and some trigonometric identities, we approximate the integral above by
\begin{equation*}
\begin{aligned}
\int_{\mathcal{A}}  \sum_{j',l}  &\frac{4\beta_{j'+l}}{a^2 N}\left[ \sin(\frac{\pi j' \xref}{a}) \cos(\frac{\pi l \xref}{a}) + \cos(\frac{\pi j' \xref}{a}) \sin (\frac{\pi l \xref}{a})\right] \sin (\frac{\pi j' x^{\rm S}
}{a}) e^{-i\frac{k^2}{\beta_{j'}} \frac{j'}{N} \frac{l}{N} \zr}\\
\times & \left[ \sin^2 (\frac{\pi j' x_1}{a}) \cos (\frac{\pi l x_1}{a}) + \cos \sin (\frac{\pi j' x_1}{a}) \sin (\frac{\pi l x_1}{a}) \right] e^{i 2\pi \sqrt{1-(\frac{j'}{N})^2} \frac{\zr - z^{\rm S}
}{\lambda} } dx_1.
\end{aligned}
\end{equation*}
Defined $\eta = z^{\rm S}  - z_{\rm r}$ and $\xi = x^{\rm S} - x_{\rm r}$. Substitution of some 
further trigonometric identities transform the above integral into
\begin{equation}
\begin{aligned}
\int_{\mathcal{A}} \frac{1}{a^2 N} \sum_{j',l} &\beta_{j'+l} 
\left\{\left(\cos(\frac{\pi j' \xi}{a}) \cos(\frac{\pi l \xref}{a}) + \sin 
(\frac{\pi j' \xi}{a}) \sin (\frac{\pi l \xref}{a}) \right) \cos (\frac{\pi 
l x_1}{a}) + R(j', \xi, l) \right\} \\
& e^{i2\pi \sqrt{1 - (\frac{j'}{N})^2} \frac{\eta}{\lambda}} 
e^{-i\frac{j'/N}{\sqrt{1-(j'/N)^2}} (\frac{k\zr}{N}) l} dx_1.
\end{aligned}
\label{eq:twoscaleps}
\end{equation}
Here, $R(j' ,\xi, l)$ consists of products of trigonometric functions, and 
for each product, one of the trigonometric functions is evaluated at $\pi 
j' \xref/a$, $\pi j' x_1/a$ or $\pi j'(\xref+x^{\rm S})/a$. Assuming that $\xref/a$ and $x_1/a$ are of 
order one, we observe that comparing with $\pi j' \xi/a$, which is 
explicitly written above and in which $\xi$ can be of order much smaller 
than one, $\pi j' \xref/a$ varies much faster as $j'$ varies. In other 
words, $\pi j' \xref/a$ can be viewed as a fast variable, and $\pi j' 
\xi/a$ with $\xi \ll a$ is a slow variable. Similarly, in the expression 
of the complex potentials and $\beta_{j' + l}$ for fixed $l$, the 
variable $j'/N$ with $N \gg 1$ involved is also a slow variable.  
Consequently, the contribution of $R(j',\xi,l)$ to the sum over $j'$ is 
negligible comparing with those of the terms explicitly written. This 
two-scale analysis also shows that the above integral is much smaller for 
large $\xi$ comparing with $\xi \ll a$.

We consider the regime $N \gg 1$, and use the continuum approximation to 
rewrite the sum over $j'$ as an integral with respect to the variable $t 
= j'/N$. We further assume $1 \ll l \ll N$, which is equivalent to say 
$a\gg (a_2 - a_1) \gg \lambda$; the second relation validates the 
linearization \eqref{eq:betalin} and allows us replacing $\beta_{j'+l}$ 
by $\beta_{j'}$, while the first relation justifies the usage of Poisson 
summation formula
\begin{equation*}
\begin{aligned}
&\sum_{l} \cos (\frac{\pi l \xref}{a}) \cos(\frac{\pi l x_1}{a}) e^{-i \pi 
\frac{t}{\sqrt{1-t^2}} (\frac{\zr}{a}) l} = \frac{\pi}{2} \sum_{m \in \Z} 
\left[ \delta(\frac{\pi(\xref + x_1)}{a} - \frac{\pi t}{\sqrt{1-t^2}} 
\frac{\zr}{a} + 2m\pi) \right. \\
+ & \delta (\frac{\pi(\xref - x_1)}{a} - \frac{\pi t}{\sqrt{1-t^2}} 
\frac{\zr}{a} + 2m\pi) + \delta (\frac{\pi(\xref - x_1)}{a} + \frac{\pi 
t}{\sqrt{1-t^2}} \frac{\zr}{a} + 2m\pi)\\
 + &  \delta \left.(\frac{\pi(\xref + x_1)}{a} + \frac{\pi t}{\sqrt{1-t^2}} 
\frac{\zr}{a} + 2m\pi) \right]
\end{aligned}
\end{equation*}
and
\begin{equation*}
\begin{aligned}
&\sum_{l} \sin (\frac{\pi l \xref}{a}) \cos(\frac{\pi l x_1}{a}) e^{-i \pi 
\frac{t}{\sqrt{1-t^2}} (\frac{\zr}{a}) l} = \frac{\pi}{2i} \sum_{m \in 
\Z} \left[ \delta(\frac{\pi(\xref + x_1)}{a} - \frac{\pi t}{\sqrt{1-t^2}} 
\frac{\zr}{a} + 2m\pi) \right. \\
+ & \delta (\frac{\pi(\xref - x_1)}{a} - \frac{\pi t}{\sqrt{1-t^2}} 
\frac{\zr}{a} + 2m\pi) - \delta (\frac{\pi(\xref - x_1)}{a} + \frac{\pi 
t}{\sqrt{1-t^2}} \frac{\zr}{a} + 2m\pi)\\
 - &  \delta \left.(\frac{\pi(\xref + x_1)}{a} + \frac{\pi t}{\sqrt{1-t^2}} 
\frac{\zr}{a} + 2m\pi) \right]
\end{aligned}
\end{equation*}
Note that we have used the fact that $k/N = \pi/a$ again. For each fixed 
$m$, let $\delta_\alpha (x_1,\xref;\zr,m)$, $\alpha = 1, \cdots, 4$, denote 
the four Dirac distributions. Define also $\tilde{\eta} = 2\pi 
\eta/\lambda$, $\tilde{\xi} = 2\pi \xi/\lambda$; Then we have $\cos(\pi 
j' \xi/a) = \cos(\tilde{\xi} t)$; the integral above becomes
\begin{equation*}
\begin{aligned}
\frac{k\pi}{2a^2} \int_0^1 \sqrt{1-t^2} \Big[ \int_{\mathcal{A}} e^{i 
(\sqrt{1-t^2} \tilde{\eta} - t \tilde{\xi})}  & \sum_{m\in \Z} (\delta_1 
+ \delta_2)(x_1,\xref; \zr, m) \\
+ &  e^{i (\sqrt{1-t^2} \tilde{\eta} + t\tilde{\xi})}  \sum_{m\in \Z}( 
\delta_3 + \delta_4)(x_1,\xref; \zr, m) dx_1 \Big] dt.
\end{aligned}
\end{equation*}
If we extend the domain of $t$ to $(-1,1)$, the above integral simplifies 
to
\begin{equation*}
\frac{k\pi}{2a^2} \int_{-1}^1 \sqrt{1-t^2} \Big[ \int_{\mathcal{A}} e^{i 
(\sqrt{1-t^2} \tilde{\eta} + t \tilde{\xi})}  \sum_{m\in \Z} (\delta_3 + 
\delta_4)(x_1,\xref; \zr, m)\Big] dt.
\end{equation*}
Integrate over $x_1$ first. The two Dirac distributions restrict the 
range of integration of $t$ to the following intervals respectively:
\begin{equation*}
\begin{aligned}
\frac{(a_1 + 2ma) - \xref}{\zr} &\le &\frac{t}{\sqrt{1-t^2}} & &\le 
\frac{(a_2 + 2ma) - \xref}{\zr},\\
\frac{(- a_2 + 2ma) - \xref}{\zr} &\le &\frac{t}{\sqrt{1-t^2}} & &\le 
\frac{(- a_1 + 2ma) - \xref}{\zr}.
\end{aligned}
\end{equation*}
To interpret these conditions, imagine that the boundaries of the 
waveguide are two mirrors, then the array $(a_1,a_2)$ has an image 
$(-a_2,-a_1)$ in the lower mirror; the two mirrors then generate a 
replica of images $(a_1 + 2ma, a_2 + 2ma)$ and $(-a_2 + 2ma, -a_1 + 
2ma)$. If we define an angle $\theta$ by
\begin{equation*}
\theta = \arctan \frac{t}{\sqrt{1-t^2}}, \quad t \in (-1,1),
\end{equation*}
then the above restrictions of the Dirac measures can be restated as
\begin{equation*}
\begin{aligned}
\arctan \frac{\xref  + a_1 + 2ma}{\zr} &\le &\theta  & &\le \arctan 
\frac{\xref + a_2 + 2ma}{\zr},\\
\arctan \frac{\xref  - a_2 + 2ma}{\zr} &\le &\theta & &\le \arctan 
\frac{\xref - a_1 + 2ma}{\zr}.
\end{aligned}
\end{equation*}
The first one restrict the angle $\theta$ to those formed by the 
reflector and the array $(a_1, a_2)$ and the images of this array. The 
second set restrict the angle $\theta$ to those formed by the reflector 
and the image array $(-a_2, -a_1)$ and the replicas.

To analyze the resulting integral, we consider the simplest set-up where: 
$(a_1, a_2)$ is centered in the cross section; the reflector is also 
centered in the cross range direction, {\itshape i.e.} $\xref = a/2$. Further we 
assume the large distance regime
$\zr \gg a$.

In such a setting, with the notation $a_c = (a_1 + a_2)/2$ and $w = (a_2 
- a_1)/2$, the integral above boils down to
\begin{equation*}
\frac{k}{2a} \sum_{m \in \Z} \int_{\arctan (ma + a_c - w - 
\xref)/\zr}^{\arctan (ma+a_c + w - \xref)/\zr} \cos^2 \theta 
e^{i(\tilde{\eta} \cos \theta + \tilde{\xi} \sin \theta)} d\theta.
\end{equation*}
Since $w \ll \zr$, for each fixed $m$, the integral is over a very small 
angle section. Hence we can approximate the integral by the value at mean 
angle times the length of the angle section. 
The mean angle $\theta_m$ in the angle section is 
$\arctan ma/\zr$. We further check that $\cos \theta_m = 
1/\sqrt{1+(ma/\zr)^2}$. 
Consequently, with $a/\zr$ set to $\triangle x$, the sum becomes
\begin{equation*}
\begin{aligned}
& &&\frac{k}{2a} \frac{2w}{\zr} \frac{\zr}{a} \sum_{m=-\infty}^\infty 
\frac{\exp\{i ( \tilde{\eta} + m \triangle x \tilde{\xi}) (\sqrt{1+(m\triangle 
x)^2})^{-1}\}}{(1+(m\triangle x)^2)^2} \triangle x\\
&\approx &&\frac{k}{2a} 
\frac{a_2 - a_1}{a} \int_{-\infty}^\infty 
\frac{e^{i\tilde{\eta}/\sqrt{1+x^2} + i\tilde{\xi}x/\sqrt{1+x^2}}}{(1+x^2)^2} dx = \frac{k}{2a} \frac{a_2 - a_1}{a} \int_{-\frac{\pi}{2}}^{\frac{\pi}{2}} \cos^2 \theta 
e^{i\tilde{\eta}\cos \theta+ i \tilde{\xi}\sin \theta} d\theta.
\end{aligned}
\end{equation*}

Finally, let us verify that the primary-primary and secondary-primary 
cross correlations do not have significant contributions to the imaging 
functional $\mathcal{I}_{{\rm LA}+}$. Let $\mathcal{I}_{+{\rm pp}}$ and 
$\mathcal{I}_{+{\rm sp}}$ denote these two terms respectively. Similar to 
\eqref{eq:si1ps}, the expectation of $\mathcal{I}_{+{\rm sp}}$ converges, as 
$\eps \to 0$, to
\begin{equation}
-\frac{\omega_0^2 \sigref \Phi_{-1}(x_{\rm s})}{8} \Big[ \frac{1}{N}\sum_{j,j'=1}^N 
\beta_j \phi_j(\xref) \phi_{j'}(x^{\rm S}
) e^{-i(\beta_j \zr + 
\beta_{j'} z^{\rm S}
)} \int_{\mathcal{A}} \phi_j(x_1) \phi_{j'}(x_1)dx_1 
\Big]^2.
\label{eq:si1sp}
\end{equation}
This function has the same form of \eqref{eq:si1ps} and can be analyzed 
in the same way. The key of these two functions is that the phase 
function in $\mathcal{I}_{+{\rm sp}}$ is a sum. As a result, the variable 
$\eta$ in \eqref{eq:twoscaleps} cannot be defined and have to be replaced 
by $\zr + z^{\rm S}
$ which is of order $a$. This renders $e^{i 2\pi 
\sqrt{1-(j'/N)^2} (\zr + z^{\rm S}
)/\lambda}$ fast varying no matter how close 
$z^{\rm S}$ is to $\zr$. Due to the averaging of fast oscillations, there is no 
significant contribution from $\mathcal{I}_{+{\rm sp}}$.

For the primary-primary component, the analog to \eqref{eq:si1ps} reads
\begin{equation}
\E \big[ \mathcal{I}_{+{\rm pp}} \big] \longrightarrow \frac{-i\Phi_{-1}(x_{\rm s})}{2a} \sum_{j} \beta_j^3 
\Big[ \int_{\mathcal{A}} \phi_j(x_1) \frac{1}{N} \sum_{j'} \sin 
(\frac{2\pi x_1}{\lambda} \frac{j}{N}) e^{-2\pi \sqrt{1-(\frac{j'}{N})^2} 
\frac{z^{\rm S}
}{\lambda}} dx_1 \Big]^2.
\end{equation}
Again, in the sum over $j'$, only fast variables are involved. In the regime $N \gg 1$, the contribution of the function above is negligible.

The term $\mathcal{I}_{{\rm LA}-}$ can be analyzed similarly. Combining the main contributions in $\mathcal{I}_{{\rm LA}+}$ and $\mathcal{I}_{{\rm LA}-}$, we obtain the desired result.
\end{proof}

%%%%%%%%%%%%%%%
\subsection{Imaging with Broadband Sources}

As we will see in Section \ref{sec:stabi}, the imaging functionals are not statistical stable if the source is time-harmonic. Hence it is required to consider
a broadband source (\ref{broadbandsource1}). 
We show that the results obtained above for the means of imaging functionals apply to the broadband setting as well.

Using Proposition \ref{prop:momentT}, the main contribution of the two moment of mode 
coupling matrix at the same frequency comes from the case when $j=m$ and 
$l=n$. Therefore,
\begin{equation*}
\begin{aligned}
\E~\eps^\alpha \Cc_{\rm pp}(\tau, x_1,x_2) \longrightarrow &
\frac{\Phi_{-1}(x_{\rm s})}{4} \sum_{j=1}^N \frac{\phi_j(x_1) \phi_j(x_2)}{\beta_j}  
\lim_{\eps \to 0} e^{-i\omega_0 \tau} \frac{1}{2\pi} \int |\hat{f}_0(h)|^2 
e^{-i\eps^\alpha \tau h} dh.
\end{aligned}
\end{equation*}

Similarly, for the primary-secondary field, we have
\begin{equation*}
\begin{aligned}
\E~\eps^{\alpha} \Cc_{\rm ps}(\tau,x_1,x_2) \longrightarrow& 
 \frac{i\omega_0^2 \Phi_{-1}(x_{\rm s}) \sigref}{8} \sum_{q,j=1}^N 
\frac{\phi_j(x_1)\phi_j(\xref) \phi_q(x_2)\phi_q(\xref)}{\beta_q \beta_j} \\
&\times \lim_{\eps\to 0} e^{i(\beta_q + \beta_j) \zr} 
e^{-i\omega_0 \tau} \frac{1}{2\pi} \int |\hat{f}_0(h)|^2  e^{i(\beta'_q + 
\beta'_j)\eps^\alpha h \zr} e^{-i\eps^\alpha h \tau}
dh.
\end{aligned}
\end{equation*}

For the secondary-primary field, we have
\begin{equation*}
\begin{aligned}
\E~\eps^{\alpha} \Cc_{\rm ps}(\tau,x_1,x_2) \longrightarrow&
\frac{-i\omega_0^2 \Phi_{-1}(x_{\rm s}) \sigref}{8} \sum_{q,j=1}^N 
\frac{\phi_q(x_1)\phi_q(\xref) \phi_j(x_2)\phi_j(\xref)}{\beta_q \beta_j} \\
&\times \lim_{\eps\to 0} e^{-i(\beta_q + \beta_j) 
\zr} e^{-i\omega_0 \tau} \frac{1}{2\pi} \int 
|\hat{f}_0(h)|^2  e^{-i(\beta'_q + 
\beta'_j)\eps^\alpha h \zr} e^{-i\eps^\alpha h \tau}
dh.
\end{aligned}
\end{equation*}

From these limits, we see that as long as $\zr$, the distance between the 
reflector and the array is much smaller than $\eps^{-\alpha}$, the 
integral in $h$ above can be approximated by the energy of the source 
(square $L^2$ norm of $f_0$). The rest parts of the limiting expectation of the cross correlation 
functions are exactly the same as the time-harmonic case. Consequently, 
the resolution analyses in the previous subsections based on the mean value of the cross-correlation migration imaging functionals remain the same.

%%%%%%%%%%%%%
%%%%%%%%%%%%%
\section{Stability Analysis of the Imaging Functionals}
\label{sec:stabi}

The key tool is the following proposition which 
analyzes the asymptotic behavior of the fourth-order moment of the transmission
coefficients in the limit 
$\eps \to 0$; see \cite[Section 20.9.3]{FGPS} or \cite[Section 8.4]{GP-SIAM07}.

\begin{proposition}
\label{prop:fourmoment}%
The expectation of four 
transmission coefficients at the same frequency has a limit as $\eps \to 
0$. In the regime $L \gg L_{\rm equip}$
we have
$$
\lim_{\eps \rightarrow 0}
\E [ \overline{T_{jl}^\eps} {T_{mn}^\eps}   
\overline{T_{j'l'}^\eps} {T_{m'n'}^\eps}   ]
\stackrel{L \gg L_{\rm equip}}{\simeq}
\left\{
\begin{array}{ll}
\frac{2}{N(N+1)}
&\mbox{ if } (j,l)=(m,n)=(j',l')=(m',n') \, ,\\
\frac{1}{N(N+1)}
&\mbox{ if }  (j,l)=(m,n) \neq (j',l')=(m',n') \, ,\\
\frac{1}{N(N+1)}
&\mbox{ if }  (j,l)=(m',n') \neq (j',l')=(m,n) \, ,\\
0 
& \mbox{ otherwise}\, .
\end{array}
\right.
$$
Let $\alpha\in (0,2)$ and $h \neq 0$.
The expectation of four 
transmission coefficients at two frequencies $\omega$ and $\omega+\eps^\alpha h$ 
has a limit as $\eps \to  0$. In the regime $L \gg L_{\rm equip}$
we have
$$
\lim_{\eps \rightarrow 0}
\E [ \overline{T_{jl}^\eps} {T_{mn}^\eps}   (\omega)
\overline{T_{j'l'}^\eps} {T_{m'n'}^\eps} (\omega+\eps^\alpha h)  ]
\stackrel{L \gg L_{\rm equip}}{\simeq}
\left\{
\begin{array}{ll}
\frac{1}{N^2}
&\mbox{ if } (j,l)=(m,n) \mbox{ and } (j',l')=(m',n') \, ,\\
0 
& \mbox{ otherwise}\, .
\end{array}
\right.
$$
\end{proposition}

The previous section shows that the mean of the imaging functional
has a peak centered at the reflector location.
The width of the peak is of the order of the wavelength.
However, the imaging functional will give the reflector location
only if it is statistically stable, that is to say, if the standard deviation of the fluctuations
of the imaging functional is smaller than the mean amplitude of the peak.

\subsection{Time-harmonic case}
We address the full aperture case in which the 
imaging functional is defined by (\ref{eq:IKMdef}).
By using (\ref{eq:espIFAp}-\ref{eq:espIFAm})
the mean of the imaging functional is 
\begin{equation}
\label{eq:meanif}
\E \big[ \mathcal{I}_{\rm FA} (x^{\rm S}, z^{\rm S}) \big] = \frac{\omega_0^2 \sigref}{4} \Phi_{-1}(\xs) \Big\{
\Re e \big( \Psi(x^{\rm S},z^{\rm S};\xref,\zr) ^2 \big) 
+ O\big( \frac{1}{kz^{\rm S}} \big)  + 
O \big( \frac{1}{N \omega_0^2 \sigref} \big) 
\Big\}  ,
\end{equation}
where $\Phi_j$ and $\Psi$ are defined by (\ref{def:Phij}) and (\ref{def:Psi}).
Here:\\
- The term with the real part comes from the contributions  of the cross correlation 
of secondary (reflected) and primary waves ${\mathcal C}_{\rm ps}$ and ${\mathcal C}_{\rm sp}$
that contain $\Psi(x^{\rm S},z^{\rm S};\xref,\zr)$ (second term in (\ref{eq:espIFAp}))
or $\Psi(x^{\rm S},-z^{\rm S};\xref,-\zr)$  (third term in (\ref{eq:espIFAm})).\\
- The term $O( 1/(k z^{\rm S}))$ comes from the contributions of the cross correlation 
of secondary (reflected) and primary waves ${\mathcal C}_{\rm ps}$ and ${\mathcal C}_{\rm sp}$
that contain $\Psi(x^{\rm S},-z^{\rm S};\xref,\zr)$ or $\Psi(x^{\rm S},z^{\rm S};\xref,-\zr)$ 
(third term in (\ref{eq:espIFAp}) and second term in (\ref{eq:espIFAm})).
In such a case there are at least a product of two these terms, which gives the $1/(kz^{\rm S})$ decay.\\
- The term $O(1/(N \omega_0^2 \sigref))$ comes the contributions of the cross correlation
of primary waves ${\mathcal C}_{\rm pp}$ (the first terms in (\ref{eq:espIFAp}) and in (\ref{eq:espIFAm})). 

The expression (\ref{eq:meanif}) is valid provided $N$ is large enough so that $N \omega_0^2 \sigref \gg 1$. 
Then it is true that the mean imaging functional is dominated by the first term in the right hand side,
which is a peak centered at the reflector location.
The mean amplitude of the peak at the reflector location is 
$$
P_{\rm peak} =  \frac{\omega_0^2 \sigref}{4} \Phi_{-1}(\xs) \Phi_0(\xref)^2  .
$$
In the continuum approximation $N\gg 1$, we have (\ref{asy:Phi-1}) and 
\begin{equation}
\label{asy:Phi0}
 \Phi_{0}(x_{\rm r}) = \frac{2}{a N} \sum_{j=1}^{N}   
  \sin^2(\frac{2\pi x_{\rm r}}{\lambda} 
\frac{j}{N})  \stackrel{N \gg 1}{\simeq} \frac{1}{a}  ,
\end{equation}
and therefore
\begin{equation}
P_{\rm peak} =  \frac{\pi \omega_0^2 \sigref}{8 a^3} .
\end{equation}

The second moment of the imaging functional can be computed  
using Proposition \ref{prop:fourmoment}
in the regime $\eps \ll 1$ and 
$L \gg L_{\rm equip}$.
\begin{eqnarray}
\nonumber
\E \big[ |\mathcal{I}_{\rm FA}  (x^{\rm S}, z^{\rm S}) |^2 \big]
&=&
\frac{N}{N+1}\big|  \E \big[ \mathcal{I}_{\rm FA}   (x^{\rm S},z^{\rm S})\big] \big|^2\\ 
\nonumber
&&+
\frac{\omega_0^4 \sigref^2 N}{32 (N+1)}\Phi_{-1}^2(\xs)
\Big\{ \Phi_1(x^{\rm S}) \Phi_{-1}(\xref) | \Psi(x^{\rm S},z^{\rm S};\xref,\zr)|^2 \\
 &&\quad +
 \Re e \big(   \Psi(x^{\rm S},z^{\rm S};\xref,\zr)^4 \big)  + O\big( \frac{1}{kz^{\rm S}} \big) + 
O \big( \frac{1}{N^2 \omega_0^4 \sigref^2} \big) 
  \Big\}  .
  \label{eq:varth1}
\end{eqnarray}
Here:\\
- The term $O( 1/(k z^{\rm S}))$ comes from the contributions of the cross correlation 
of secondary (reflected) and primary waves ${\mathcal C}_{\rm ps}$ and ${\mathcal C}_{\rm sp}$
that contain $\Psi(x^{\rm S},-z^{\rm S};\xref,\zr)$ or $\Psi(x^{\rm S},z^{\rm S};\xref,-\zr)$ 
(in such a case there are at least a product of two these terms, which gives the $1/(kz^{\rm S})$ decay).\\
- The term $O(1/(N \omega_0^2 \sigref)^2)$ comes from the contribution of the cross correlation of primary waves ${\mathcal C}_{\rm pp}$ that can be computed in a more quantitative way:
$$
O \big( \frac{1}{N^2 \omega_0^4 \sigref^2} \big) 
= \frac{4}{N^2 \omega_0^4 \sigref^2} \Phi_1(x^{\rm S})^2
+ o \big( \frac{1}{N^2 \omega_0^4 \sigref^2} \big)   .
$$

The variance of the imaging functional at the reflector location is therefore
\begin{eqnarray*}
{\rm Var} \big( \mathcal{I}_{\rm FA} (\xref,  \zr)  \big) = P_{\rm peak}^2   \Big\{  \frac{1}{2} + \frac{1}{2} 
\frac{\Phi_1(\xref) \Phi_{-1}(\xref)}{\Phi_0^2(\xref)}   + O\big( \frac{1}{kz^{\rm S}
} \big) + 
O \big( \frac{1}{N^2 \omega_0^4 \sigref^2} \big) 
  \Big\}  .
\end{eqnarray*}
In the continuum approximation $N\gg 1$, we have
(\ref{asy:Phi-1}), (\ref{asy:Phi0}), and 
\begin{equation}
\label{asy:Phi1}
 \Phi_{1}(x_{\rm r}) = \frac{2}{a N} \sum_{j=1}^{N} \beta_j
  \sin^2(\frac{2\pi x_{\rm r}}{\lambda} 
\frac{j}{N})  \stackrel{N \gg 1}{\simeq} \frac{1}{a} \int_0^1 \sqrt{1-x^2} dx = \frac{\pi}{4a},
\end{equation}
and therefore
\begin{equation}
{\rm Var} \big( \mathcal{I}_{\rm FA} (\xref, \zr)  \big) = P_{\rm peak}^2  \Big\{ \frac{1}{2} + \frac{1}{2} \frac{\pi^2}{8}  \Big\}  .
  \label{eq:varth2}
\end{equation}

To summarize:\\
1) The typical amplitude of the fluctuations of the imaging  functional for $(x^{\rm S},z^{\rm S}) = (\xref,\zr)$
({\itshape i.e.} at the reflector location) is $P_{\rm peak}$ (as shown by (\ref{eq:varth2})).\\
2) The typical amplitude of the fluctuations of the imaging  functional for $|(x^{\rm S},z^{\rm S})-(\xref,\zr)| \gg \lambda_0$
({\itshape i.e.} away from the reflector location)
is $P_{\rm peak} \big( (\lambda_0/|z^{\rm S}
-\zr|)^{1/2} \wedge ( \lambda_0 |x^{\rm S}
-\xref|)^{3/2} + \lambda_0^2/(N \sigref) \big)$ (as shown by (\ref{eq:varth1})).

The second result shows that the fluctuations of the image far from the main peak location are of the  order of $ \lambda_0^2/(N \sigref)$ relatively to the amplitude of the main peak. They are due to the contributions of the primary cross correlation. Provided the number of modes is large enough $N \gg \lambda_0^2 / \sigref $,
they are small.

The first result shows that the amplitude of the peak at the reflector location has relative fluctuations of order one. This is 
due to the fact that the reflector is illuminated by a field whose amplitude is randomly spatially varying, so that the reflected energy is proportional to the squared amplitude of the primary field at the reflector location, which is a random quantity. This is the origin of the statistical instability in the time-harmonic case.\\

\subsection{Broadband case}
We know that the frequency coherence radius $\Omega_c$
in a waveguide with length $L/\eps^2$
is of the order of $\eps^2$ (see \cite[Proposition 20.7]{FGPS} or \cite[Proposition 6.3]{GP-SIAM07}).
As a result, as soon as a broadband source with a bandwidth larger than $\Omega_c$
is used, then the field is the superposition of decorrelated frequency components.
As a consequence the field is self-averaging in the time domain.

More exactly, from the expressions of the cross correlations in the broadband case given in Subsection \ref{subsec:crossbroad}, the mean 
and the variance of the imaging functional are of the form
\begin{eqnarray*}
\E \big[  \mathcal{I}_{\rm FA} (x^{\rm S},  z^{\rm S}) \big] &=&
\int dh  |\hat{f}_0(h)|^2  
\sum_{j,l,m,n}
 \E†\big[ \overline{T_{jl}^\eps} T^\eps_{mn}(\omega_0+\eps^\alpha h)\big] 
c_{j,l,m,n} (x^{\rm S},  z^{\rm S}) , \\
{\rm Var} \big(  \mathcal{I}_{\rm FA} (x^{\rm S},  z^{\rm S}) \big) &=&
\iint dh dh' |\hat{f}_0(h)|^2 |\hat{f}_0(h')|^2\sum_{j,l,m,n,j',l',m',n'} \\
&&
\Big\{ \E†\big[ \overline{T_{jl}^\eps} T^\eps_{mn}(\omega_0+\eps^\alpha h)
T_{j'l'}^\eps \overline{T^\eps_{m'n'}}(\omega_0+\eps^\alpha h') \big]   \\
&& \hspace*{-0.9in} - 
 \E†\big[ \overline{T_{jl}^\eps} T^\eps_{mn}(\omega_0+\eps^\alpha h)\big]
\E \big[ T_{j'l'}^\eps \overline{T^\eps_{m'n'}}(\omega_0+\eps^\alpha h') \big] \Big\}
c_{j,l,m,n} (x^{\rm S},  z^{\rm S}) \overline{c_{j',l',m',n'}} (x^{\rm S},  z^{\rm S}) ,
\end{eqnarray*}
where $c_{j,l,m,n} (x^{\rm S},  z^{\rm S})$ is a shorthand for the deterministic coefficient
that contains the phase and mode amplitudes.

First, since $ \E†\big[ \overline{T_{jl}^\eps} T^\eps_{mn}(\omega_0+\eps^\alpha h)\big] $
is independent on $h$ to leading order (because $\alpha>0$), the mean  satisfies
$$
\E \big[ \mathcal{I}_{\rm FA} (x^{\rm S},  z^{\rm S})  \big]_{\rm broadband} 
\sim
\E \big[ \mathcal{I}_{\rm FA} (x^{\rm S},  z^{\rm S}) \big]_{\rm narrowband}  ,
$$
as already noticed.

Second,
the term in the curly brackets in the expression of the variance is vanishing if $\eps^\alpha |h-h' |$ is larger than the  frequency coherence radius $\Omega_c$. 
So the double integral in $(h,h')$ is reduced to a domain that has the form
of a thin diagonal band, whose thickness is limited by the frequency coherence radius $\Omega_c$.
As a result we obtain that
\begin{eqnarray*}
{\rm Var} \big(  \mathcal{I}_{\rm FA} (x^{\rm S},  z^{\rm S}) \big) &\sim&
\iint_{\eps^\alpha |h-h'|  \leq \Omega_c} dh dh' |\hat{f}_0(h)|^2 |\hat{f}_0(h')|^2
\sum_{j,l,m,n,j',l',m',n'} \\
&&
\Big\{ \E†\big[ \overline{T_{jl}^\eps} T^\eps_{mn}(\omega_0)
T_{j'l'}^\eps \overline{T^\eps_{m'n'}}(\omega_0) \big]   \\
&& \hspace*{-0.9in} - 
 \E†\big[ \overline{T_{jl}^\eps} T^\eps_{mn}(\omega_0)\big]
\E \big[ T_{j'l'}^\eps \overline{T^\eps_{m'n'}}(\omega_0) \big] \Big\}
c_{j,l,m,n} (x^{\rm S},  z^{\rm S}) \overline{c_{j',l',m',n'}} (x^{\rm S},  z^{\rm S}) ,
\end{eqnarray*}
or more simply
\begin{eqnarray*}
{\rm Var} \big( \mathcal{I}_{\rm FA} (x^{\rm S},  z^{\rm S})  \big)_{\rm broadband} 
\sim
{\rm Var} \big( \mathcal{I}_{\rm FA} (x^{\rm S},  z^{\rm S}) \big)_{\rm narrowband}
\frac{\Omega_c}{B}  ,
\end{eqnarray*}
where $B \sim \eps^\alpha$ is the bandwidth of the source that is larger than the 
frequency coherence radius $\Omega_c \sim \eps^2$.
We had seen that the use of broadband sources does not affect the resolution of the imaging
functional but it ensures its statistical stability.
Provided the bandwidth is larger than the  frequency coherence radius, the typical amplitude of the 
fluctuations of the imaging functional is smaller than the amplitude $P_{\rm peak}$ of the main peak at the reflector location, 
and therefore the reflector can be localized.

%%%%%%%%%%%%%%
%%%%%%%%%%%%%%
\section{Conclusions}
\label{sec:discu}
In this paper we have shown that migration of the cross correlations
of the data recorded by a passive receiver array 
can allow for diffraction-limited imaging of the reflector 
 in a random waveguide even though the sources
are very far from the reflector, provided the receivers are close enough from it.
% the data recorded at a passive receiver array  can be used in order to  improve migration techniques in a random waveguide.
The statistical stability of the imaging functional is ensured by the use of broadband sources.
The resolution properties are ensured by the waveguide geometry:
even when the receiver array does not span the whole cross section of the waveguide, the
width of the point spread function of the imaging functional is of the order of the wavelength,
provided the diameter of the array is larger than the wavelength.

This paper has addressed the case of a two-dimensional waveguide with Dirichlet boundary conditions,
but the conclusions should be qualitatively the same for fairly general situations, 
when addressing three-dimensional waveguides,
with Neumann, Dirichlet or mixed boundary conditions, with random fluctuations of the index of refraction or of
 the boundaries  as in \cite{alonso,gomez}.

\section*{Acknowledgements}
%[START OF REVISION]
The authors would like to thank the anonymous referees for their careful reading of the manuscript and helpful comments. This work was supported  by ERC Advanced Grant Project MULTIMOD-267184.
%[END OF REVISION]

\end{document}